\documentclass[leqno,11pt]{article}

\usepackage{amssymb}
\usepackage{amsmath}
\usepackage{amsthm}
\usepackage[dvipdfmx]{graphicx}
\usepackage{bm}
\usepackage{cases}
\usepackage{mathrsfs}

\makeatletter
    
    \@addtoreset{equation}{section}
\makeatother

 \newtheorem{thm}{Theorem}[section]
 
 \newtheorem{lem}[thm]{Lemma}
 \newtheorem{prop}[thm]{Proposition}
 
 \newtheorem{assu}[thm]{Assumption}

\begin{document}

\begin{center}
\bf\Large
A Mathematical Justification of a Thin Film Approximation \\
for the Flow down an Inclined Plane 
\end{center}

\begin{center}
Hiroki Ueno and Tatsuo Iguchi
\end{center}

\begin{abstract}
We consider a two-dimensional motion of a thin film flowing down an inclined plane 
under the influence of the gravity and the surface tension. 
In order to investigate the stability of such flow, we often use a thin film approximation, which is an approximation obtained by the perturbation expansion with respect to the aspect ratio of the film. The famous example of the approximate equations are the Burgers equation, Kuramoto--Sivashinsky equation, KdV--Burgers equation, KdV--Kuramoto--Sivashinsky equation, and so on. 
In this paper, we give a mathematically rigorous justification of a thin film approximation by establishing an error estimate between the solution of the Navier--Stokes equations and  those of approximate equations. 
\end{abstract}

\section{Introduction}
In this paper, we consider a two-dimensional motion of liquid film of a viscous and incompressible fluid flowing down an inclined plane under the influence of the gravity and the surface tension on the interface. 
The motion can be mathematically formulated as a free boundary problem for the incompressible Navier--Stokes equations. 
We assume that the domain $\Omega(t)$ occupied by the liquid at time $t\ge0$, the liquid surface $\Gamma(t)$, 
and the rigid plane $\Sigma$ are of the forms 
$$
\left\{
\begin{array}{l} 
\Omega (t)=\{(x,y)\in\mathbb{R}^2\ |\ 0<y<h_0+\eta(x,t)\}, \\
\Gamma(t)=\{(x,y)\in\mathbb{R}^2\ |\ y=h_0+\eta(x,t)\}, \\
\Sigma=\{(x,y)\in\mathbb{R}^2|\ y=0\},
\end{array}
\right.
$$
where $h_0$ is the mean thickness of the liquid film and $\eta(x,t)$ is the amplitude of the liquid surface. 
Here we choose a coordinate system $(x, y)$ so that $x$ axis is pointed to the streamwise direction  and $y$ axis is normal to the plane. 
We consider fluctuations of the Nusselt flat film solution, which is the stationary laminar flow given by 
\begin{equation}\label{laminar flow}
\eta_1=0, \quad u_1=(\rho g\sin\alpha/2\mu)(2h_0y-y^2), \quad v_1=0, \quad p_1=p_0-\rho g\cos\alpha(y-h_0),
\end{equation}
where $\rho$ is a constant density of the liquid, $g$ is the acceleration of the gravity, 
$\alpha$ is the angle of inclination, $\mu$ is the shear viscosity coefficient, and $p_0$ is an  atmospheric pressure. 
Throughout this paper, we assume that the flow is $l_0$-periodic in the streamwise direction $x$.
Rescaling the independent and dependent variables by using $h_0$, $l_0$, the typical amplitude of the liquid surface $a_0$, $U_0=\rho gh_0^2\sin\alpha/2\mu$, and $P_0=\rho gh_0\sin\alpha$, the equations are written in the non-dimensional form
\begin{equation}\label{ns}
\left\{
 \begin{array}{lcl}
  \delta\bm{u}_t+\bigl((\bar{\bm{u}}+\varepsilon\bm{u})\cdot\nabla_\delta\bigr)\bm{u}
   +(\bm{u}\cdot\nabla_\delta)\bar{\bm{u}}+\dfrac{2}{{\rm R}}\nabla_\delta p
   -\dfrac{1}{{\rm R}}\Delta_\delta\bm{u}=\bm{0}
  & \mbox{in} & \Omega_{\varepsilon}(t), \; t>0, \\ 
  \nabla_\delta\cdot\bm{u} = 0 
  & \mbox{in} & \Omega_{\varepsilon}(t), \; t>0,
 \end{array}
\right.
\end{equation}
\begin{equation}\label{bc1}
\left\{
 \begin{array}{lcl}
  \bigl(\bm{D}_\delta(\varepsilon\bm{u}+\bar{\bm{u}})-\varepsilon p\bm{I}\bigr)\bm{n} & & \\
  \quad
   = \bigg(-\dfrac{1}{\tan\alpha}\varepsilon\eta
    +\dfrac{\delta^2{\rm W}}{\sin\alpha}\dfrac{\varepsilon\eta_{xx}}{(1+(\varepsilon\delta\eta_x)^2)^{\frac32}}\bigg)
    \bm{n}
    & \mbox{on} & \Gamma_{\varepsilon}(t), \; t>0, \\[3mm]
  \eta_t+\bigl(1-(\varepsilon\eta)^2+\varepsilon u\bigr)\eta_x-v = 0 
    & \mbox{on} & \Gamma_{\varepsilon}(t), \; t>0, 
 \end{array}
\right.
\end{equation}
\begin{equation}\label{bc0}
\bm{u}=\bm{0} \quad\mbox{on}\quad \Sigma, \; t>0.
\end{equation}
Here,  $\delta, \varepsilon, {\rm R},$ and ${\rm W}$ are non-dimensional parameters defined by
$$
\delta=\frac{h_0}{l_0},\quad\varepsilon=\frac{a_0}{h_0},\quad {\rm R}=\frac{\rho U_0 h_0}{\mu}, \quad {\rm W}=\frac{\sigma}{\rho g h_0^2},
$$
where $\sigma$ is the surface tension coefficient.
Note that $\delta$ is the aspect ratio of the film, $\varepsilon$ represents the magnitude of nonlinearity, ${\rm R}$ is the Reynolds number, and ${\rm W}$ is the Weber number.  
Moreover, we used notations $\bm{u}=(u,\delta v)^{\rm T}, \bar{\bm{u}}=(\bar{u},0)^{\rm T}, \bar{u}=2y-y^2, \nabla_\delta=(\delta\partial_x, \partial_y)^{\rm T}, \Delta_\delta=\nabla_\delta\cdot\nabla_\delta, \bm{D}_\delta\bm{f}=\frac{1}{2}\bigl\{\nabla_\delta(\bm{f}^{\rm T})+\bigl(\nabla_\delta(\bm{f}^{\rm T})\bigr)^{\rm T}\bigr\}$, and $\bm{n}=(-\varepsilon\delta\eta_x,1)^{\rm T}$.
In this scaling, the liquid domain $\Omega_{\varepsilon}(t)$ and the liquid surface $\Gamma_\varepsilon(t)$ are of the forms 
$$
\left\{
\begin{array}{l} 
\Omega_{\varepsilon}(t)=\{(x,y)\in\mathbb{R}^2\ |\ 0<y<1+\varepsilon\eta(x,t)\}, \\
\Gamma_{\varepsilon}(t)=\{(x,y)\in\mathbb{R}^2\ |\ y=1+\varepsilon\eta(x,t)\}.
\end{array}
\right.
$$

\noindent
Concerning a mathematical analysis of the problem in the case of $\delta=\varepsilon=1$, Teramoto \cite{Teramoto} showed that the initial value problem to the Navier--Stokes equations \eqref{ns}--\eqref{bc0} has a unique solution globally in time under the assumptions that the Reynolds number and the initial data are sufficiently small. 
Nishida, Teramoto, and Win \cite{Nishida} showed the exponential stability of the Nusselt flat film solution under the assumptions that the angle of inclination is sufficiently small and $x\in\mathbb{T}$ in addition to the assumptions in \cite{Teramoto}. 
Furthermore, Uecker \cite{Uecker} studied the asymptotic behavior for $t\to\infty$ of the solution in the case of  $x\in\mathbb{R}$ and showed that the perturbations of the Nusselt flat film solution decay like the self-similar solution of the Burgers equation under the assumptions that the initial data are sufficiently small and ${\rm R}<{\rm R}_c$. Here, ${\rm R}_c=\frac45\frac1{\tan\alpha}$ is the critical Reynolds number given by Benjamin \cite{Benjamin}. 
On the other hand, Ueno, Shiraishi, and Iguchi \cite{Ueno} derived a uniform estimate for the solution of \eqref{ns}--\eqref{bc0} with respect to $\delta$ when the Reynolds number, the angle of inclination, and the initial data are sufficiently small.

Benney \cite{Benney} derived the following single nonlinear evolution equation 
\begin{align}
\eta_t&+2(1+\varepsilon\eta)^2\eta_x-\frac{8}{15}({\rm R}_c-{\rm R})\delta\eta_{xx}+C_1\delta^2\eta_{xxx}\label{Benney}\\
&+C_2\varepsilon\delta(\eta\eta_{xx}+\eta_x^2)+\frac23\frac{\rm W}{\sin\alpha}\delta^3\eta_{xxxx}
 =O(\delta^3+\varepsilon^2\delta+\varepsilon\delta^2)\notag
\end{align}
by using the method of perturbation expansion of the solution $(u, v, p)$ with respect to $\delta$ under the thin film regime $\delta\ll1$. Here, $C_1=C_1({\rm R}, \alpha)$ and $C_2=C_2({\rm R}, \alpha)$ are constants independent of $\delta, \varepsilon$, and ${\rm W}$.   Explicit forms of $C_1$ and $C_2$ will be given in  Section 3. 
Many approximate equations are obtained from \eqref{Benney} by assuming that parameters $\varepsilon, {\rm W}$, and $\rm R$ have appropriate orders in $\delta$. In the following, we assume  $\varepsilon=\delta$ and ${\rm R}<{\rm R}_c$ and set
\begin{equation}\label{etazeta}
\eta(x, t)=\zeta(x-2t, \varepsilon t). 
\end{equation}

\medskip

\noindent
{\bf{I. Burgers equation}}

Assuming ${\rm W}_1\le{\rm W}\le \delta^{-1}{\rm W}_2$ in \eqref{Benney}, we have
\[
\eta_t+2\eta_x+4\varepsilon\eta\eta_x-\frac{8}{15}({\rm R}_c-{\rm R})\delta\eta_{xx}=O(\delta^2).
\]
Plugging \eqref{etazeta} in the above equation and passing to the limit  $\varepsilon=\delta\to0$, we obtain 
\begin{equation}\label{zeta I}
\zeta_\tau+4\zeta\zeta_x-\dfrac{8}{15}({\rm{R}}_c-{\rm{R}})\zeta_{xx}=0.
\end{equation}

\noindent
{\bf II. Burgers equation with a fourth order dissipation term}

Assuming ${\rm W}=\delta^{-2}{\rm W}_2$ in \eqref{Benney},  we have 
\[
\eta_t+2\eta_x+4\varepsilon\eta\eta_x-\frac{8}{15}({\rm R}_c-{\rm R})\delta\eta_{xx}+\frac23\frac{{\rm W}_2}{\sin\alpha}\delta\eta_{xxxx}=O(\delta^2).
\]
Plugging \eqref{etazeta}  in the above equation and passing to the limit  $\varepsilon=\delta\to0$, we obtain 
\begin{equation}\label{zeta II}
\zeta_\tau+4\zeta\zeta_x-\dfrac{8}{15}({\rm{R}}_c-{\rm{R}})\zeta_{xx}+\dfrac{2}{3}\dfrac{\rm W_2}{\sin\alpha}\zeta_{xxxx}=0.
\end{equation}

\noindent
{\bf III. Burgers equation with dispersion and nonlinear terms}

Assuming ${\rm W}_1\le{\rm W}\le {\rm W}_2$ in \eqref{Benney}, we have
\[
\eta_t+2\eta_x+4\varepsilon\eta\eta_x-\frac{8}{15}({\rm R}_c-{\rm R})\delta\eta_{xx}+C_1\delta^2\eta_{xxx}+C_2\varepsilon\delta(\eta\eta_{xx}+\eta_x^2)
+2\varepsilon^2\eta^2\eta_x=O(\delta^3).
\]
Plugging \eqref{etazeta} in the above equation and neglecting the terms of $O(\delta^3)$, we obtain 
\begin{equation}\label{zeta III}
\zeta_{\tau}+4\zeta\zeta_x-\dfrac{8}{15}({\rm{R}}_c-{\rm{R}})\zeta_{xx}+\delta\big\{C_1\zeta_{xxx}+C_2\big(\zeta\zeta_{xx}+\zeta_x^2\big)+2\zeta^2\zeta_x\big\}=0.
\end{equation}

\noindent
{\bf IV. Burgers equation with fourth order dissipation, dispersion, and nonlinear terms}

Assuming ${\rm W} = \delta^{-1}{\rm W}_2$ in \eqref{Benney}, we have
\begin{align*}
&\eta_t+2\eta_x+4\varepsilon\eta\eta_x-\frac{8}{15}({\rm R}_c-{\rm R})\delta\eta_{xx}\\
&\phantom{\eta_t}+C_1\delta^2\eta_{xxx}+C_2\varepsilon\delta(\eta\eta_{xx}+\eta_x^2)+2\varepsilon^2\eta^2\eta_x+\frac23\frac{{\rm W}_2}{\sin\alpha}\delta^2\eta_{xxxx}=O(\delta^3).\notag
\end{align*}
Plugging \eqref{etazeta} in the above equation and neglecting the terms of $O(\delta^3)$, we obtain 
\begin{align}
&\zeta_{\tau}+4\zeta\zeta_x-\dfrac{8}{15}({\rm{R}}_c-{\rm{R}})\zeta_{xx}\label{zeta IV}\\
&\quad+\delta\bigg\{C_1\zeta_{xxx}
+C_2\big(\zeta\zeta_{xx}+\zeta_x^2\big)+2\zeta^2\zeta_x+\dfrac23\dfrac{{\rm W}_2}{\sin\alpha}\zeta_{xxxx}\bigg\}=0.\notag
\end{align}

We remark that \eqref{zeta III} and \eqref{zeta IV} are higher order approximate equations to Burgers equation \eqref{zeta I}. 
In this paper, we assume ${\rm R}\ll{\rm R}_c$ in order to use a uniform estimate in $\delta$ for the solution of the Navier--Stokes equations, which is a severe restriction.
Uniform estimates in $\delta$ for the solution play a most important role in the  justification for these approximation. 
Here if we could assume ${\rm R}>{\rm R}_c$, then \eqref{zeta II} would be the Kuramoto--Sivashinsky equation (see \cite{Kuramoto}, \cite{Sivashinsky1}, and \cite{Sivashinsky2}).
If we could assume ${\rm R}_c-{\rm R}=\delta\tilde{\rm R}>0$, then we would obtain the $\delta$-independent KdV--Burgers equation 
\begin{equation}\label{KdV-Burgers zeta 2}
\zeta_{\tau}+4\zeta\zeta_x-\frac{8\rm \tilde{R}}{15}\zeta_{xx}+C_1 \zeta_{xxx}=0
\end{equation}
by  plugging \eqref{etazeta} in \eqref{Benney} and passing to the limit  $\varepsilon=\delta^2\to0$  under the assumption ${\rm W}_1\le{\rm W}\le {\rm W}_2$. 
Moreover if we could assume ${\rm R}_c-{\rm R}=-\delta\tilde{\rm R}<0$,  we would obtain the $\delta$-independent KdV--Kuramoto--Sivashinsky equation  
\begin{equation}\label{KdV-KS zeta 2}
\zeta_{\tau}+4\zeta\zeta_x+\frac{8\rm \tilde{R}}{15}\zeta_{xx}+C_1\zeta_{xxx}+\frac{2}{3}\frac{\rm W_2}{\sin\alpha}\zeta_{xxxx}=0
\end{equation}
by plugging \eqref{etazeta} in \eqref{Benney} and passing to the limit $\varepsilon=\delta^2\to0$  under the assumption ${\rm W} = \delta^{-1}{\rm W}_2$. 
More details or a list of useful references about the thin film approximation can be found in 
\cite{Chang, Craster, Kalli, Oron, Ueno}.

\medskip

In this paper, we will give a mathematically rigorous justification of these thin film approximations by establishing an error estimate between the solution of the Navier--Stokes equations \eqref{ns}--\eqref{bc0} and those of the approximate equations \eqref{zeta I}--\eqref{zeta IV}.
We note that  we cannot just yet justify the Kuramoto--Sivashinsky equation, the $\delta$-independent KdV--Burgers equation  \eqref{KdV-Burgers zeta 2}, and the KdV--Kuramoto--Sivashinsky equation  \eqref{KdV-KS zeta 2} because without the assumption $\rm R\ll R_c$ we have not yet obtain a uniform estimate in $\delta$ for the solution. 
We also remark that Bresch and Noble \cite{Bresch} justified the shallow water model by proving that remainder terms converges to 0 as $\delta\to0$ (see also \cite{Bresch2}).

\medskip
The plan of this paper is as follows. 
In Section 2, we give our main theorem after we transform the problem in a time dependent domain to a problem in a time independent domain. 
In Section 3, we derive approximate solutions by using  Benney's method.
In Section 4, we recall the energy estimate for the solution of the Navier--Stokes equations obtained in \cite{Ueno}.
Finally, we give an error estimates in Section 5.

\bigskip
\noindent
{\bf Notation}. \ 
We put $\Omega=\mathbb{T}\times(0,1)$ and $\Gamma=\mathbb{T}\times\{y=1\}$, 
where $\mathbb{T}$ is the flat torus $\mathbb{T}=\mathbb{R}/\mathbb{Z}$. 
For a Banach space $X$, we denote by $\|\cdot\|_X$ the norm in $X$. 
For $1\le p\le \infty$, we put $\|u\|_{L^p}=\|u\|_{L^p(\Omega)}$, $\|u\|=\|u\|_{L^2}$, 
$|u|_{L^p}=\|u(\cdot, 1)\|_{L^p(\mathbb{T})}$, and $|u|_0=|u|_{L^2}$. 
We denote by $(\cdot, \cdot)_\Omega$ and $(\cdot, \cdot)_\Gamma$ the inner products of 
$L^2(\Omega)$ and $L^2(\Gamma)$, respectively. 
For $s\ge 0$, we denote by $H^s(\Omega)$ and $H^s(\Gamma)$ the $L^2$ Sobolev spaces of order $s$ 
on $\Omega$ and $\Gamma$, respectively. 
The norms of these spaces are denoted by $\|\cdot\|_s$ and $|\cdot|_s$. 
For a function $u=u(x,y)$ on $\Omega$, 
a Fourier multiplier $P(D_x)$ ($D_x=-{\rm{i}}\partial_x$) is defined by $(P(D_x)u)(x,y)=
\sum_{n\in\mathbb{Z}}P(n)\hat{u}_n(y){\rm{e}}^{2\pi {\rm{i}}nx}$,
where  $\hat{u}_n(y)=\int_0^1 u(x,y){\rm{e}}^{-2\pi {\rm i}nx}\,{\rm{d}}x$ is the Fourier coefficient in $x$. 
We put $\partial_y^{-1} f(x,y)=-\int_y^1 f(x,z){\rm{d}} z$ and $D_\delta^kf=\{(\delta\partial_x)^i\partial_y^j f \,|\, i+j=k\}$. 
$f\lesssim g$ means that there exists a non-essential positive constant $C$ such that $f\le Cg$ holds.

\section{Main results}
We rewrite the system \eqref{ns}--\eqref{bc0} according to \cite{Beale, Ueno}. Transforming the problem in the moving domain $\Omega(t)$ to a problem in the fixed domain $\Omega$ by using an appropriate diffeomorphism, and  introducing new unknown function $(u', v', p')$ to keep the solenoidal condition, we obtain 
\begin{equation}\label{NS-2}
\left\{
 \begin{array}{lcl}
\delta(u_t+\bar{u}u_x+\bar{u}_yv)+\dfrac2{\rm R}\delta p_x-\dfrac1{\rm R}(\delta^2 u_{xx}+u_{yy})=\delta^2 f_1& \mbox{in} & \Omega, \; t>0, \\[3mm]
\delta^2(v_t+\bar{u}v_x)+\dfrac2{\rm R} p_y-\dfrac1{\rm R}\delta(\delta^2 v_{xx}+v_{yy})=\delta^2 f_2& \mbox{in} & \Omega, \; t>0, \\[3mm]
  u_x+v_y=0 & \mbox{in} & \Omega, \; t>0,
 \end{array}
\right.
\end{equation}
\begin{equation}\label{BC1-2}
\left\{
 \begin{array}{lcl}
  \delta^2v_x+u_y-2(1+\varepsilon\eta)^2\eta=\delta^3h_1 & \mbox{on} & \Gamma, \; t>0, \\[1mm]
  p-\delta v_y-\dfrac{1}{\tan\alpha}\eta+\dfrac{\delta^2{\rm W}}{\sin\alpha}\eta_{xx}=\delta^2h_2
   & \mbox{on} & \Gamma, \; t>0, \\[2mm]
\eta_t+\eta_x-v=\delta^2\eta^2\eta_x=:\delta^2h_3 & \mbox{on} & \Gamma, \; t>0,
\end{array}
\right.
\end{equation}
\begin{equation}\label{BC0-2}
u=v=0 \quad\mbox{on}\quad \Sigma, \; t>0,
\end{equation}
where we dropped the prime sign in the notation and  $f_1, f_2,  h_1$, and $h_2$ are collections of nonlinear terms. See \cite{Ueno} for more details on the explicit form of these nonlinear terms. 
In the following, we will consider the initial value problem to \eqref{NS-2}--\eqref{BC0-2} 
under the initial conditions 
\begin{equation}\label{initial}
\eta|_{t=0}=\eta_0\quad {\rm on}\ \ \Gamma, \quad (u, v)^{\rm T}|_{t=0}=(u_0, v_0)^{\rm T}\quad{\rm in}\ \ \Omega. 
\end{equation}
Here, we assume $\int_0^1\eta_0(x){\rm d}x=0$ and denote $h_1$ determined from initial data by $h_1^{(0)}$.

We impose the following assumption on the non-dimensional parameters and initial data.
\begin{assu}\label{assumption}\ Let ${\rm R_0, R_1}, \alpha_0, {\rm W}_1, c_0$, and $M$ be positive constants and $m\ge2$ be an integer.

\begin{enumerate}
\item[(1)] Conditions for parameters

Parameters ${\rm R}, \alpha, {\rm W}, \delta$, and $\varepsilon$ satisfy
\[
{\rm R}_1\le{\rm R}\le{\rm R}_0,\quad 0<\alpha\le\alpha_0, \quad {\rm W}_1\le{\rm W}, \quad 0<\varepsilon=\delta\le1.
\]

\item[(2)] Smallness of initial data

Initial data $(\eta_0, u_0, v_0)$ and parameters {\rm W} and $\delta$ satisfy
\begin{align}
&|(1+\delta|D_x|)^2\eta_0|_2+\|(1+|D_x|)^2(u_0, \delta v_0)^{\rm T}\| +\|(1+|D_x|)^2D_{\delta}(u_0, \delta v_0)^{\rm T}\|\notag\\
&\quad+\|(1+|D_x|)^2D_{\delta}^2(u_0, \delta v_0)^{\rm T}\|+\delta^2{\rm W}|(1+\delta|D_x|)\eta_{0x}|_3+\sqrt{\delta^2{\rm W}}\|(1+|D_x|)^2\delta v_{0xy}\|\le c_0.\notag
\end{align}
\item[(3)] Regularity of  initial data

Initial data $(\eta_0, u_0, v_0)$ satisfies
\[
\|(1+|D_x|)^{m+1}(u_0, v_0)^{\rm T}\|_{H^2(\Omega)}+|\eta_0|_{m+4}\le M.
\]
\item[(4)] Compatibility conditions

Initial data $(\eta_0, u_0, v_0)$  and  parameters $\delta$ and $\varepsilon$ satisfy
\[
\left\{
\begin{array}{lcl}
 u_{0x}+v_{0y}=0 & \mbox{in} & \Omega, \\
 u_{0y}+\delta^2 v_{0x}-2(1+\varepsilon\eta_0)^2\eta_0=\delta^3h_1^{(0)} & \mbox{on} & \Gamma, \\
 u_0=v_0=0 & \mbox{on} & \Sigma.
\end{array}
\right.
\]
\end{enumerate}
\end{assu}
\noindent
{\bf{Remark 2.1.}}
Under the assumption that there exist small positive constants ${\rm R}_0, \alpha_0, $ and $c_0$ such that Assumption \ref{assumption} is fulfilled, Ueno, Shiraishi, and Iguchi \cite{Ueno} proved the global in time uniform estimate with respect to $\delta$ for the solution of the Navier--Stokes equations \eqref{NS-2}--\eqref{initial}. See also Proposition \ref{uniform estimate} in this paper.

\bigskip

For later use, we define the norm of a difference  between the solution $(\eta^\delta, u^\delta, v^\delta, p^\delta)$ of the Navier--Stokes equations \eqref{NS-2}--\eqref{initial}  and the solution  $\zeta$ of the approximate equations as 
\begin{align}\label{difference}
\mathcal{D}(t; \zeta, u, v,p):=&|\eta^\delta(t)-\zeta(\cdot-2t, \varepsilon t)|^2_0+\|(1+|D_x|)^{m}(u^\delta-u)(t)\|^2\\
&+\|(1+|D_x|)^{m-1}(v^\delta-v)(t)\|^2+\|(1+|D_x|)^{m-1}(p^\delta-p)(t)\|^2,\notag
\end{align}
\noindent
where $(u, v, p)$ is an approximate solution constructed from $\zeta$.
Let $\zeta^I, \zeta^{II}, \zeta^{III}$, and $\zeta^{IV}$ be the solution of \eqref{zeta I}--\eqref{zeta IV} under the initial condition $\zeta|_{\tau=0}=\eta_0$, respectively. 

Now we are ready to state our main results in this paper. Note that the definitions of $u^I, v^I, p^I, u^{II}, \ldots$ appeared in the following statement will be given in Section 5.

\begin{thm}\label{main theorem}
There exist small positive constants ${\rm R}_0$ and $\alpha_0$ such that the following statement holds: 
Let $m$ be an integer satisfying $m\ge2$, $0<{\rm R}_1\le {\rm R}_0$, $0<{\rm W}_1\le{\rm W}_2$, and $0<\alpha\le \alpha_0$. There exists small positive constant $c_0$ such that
if the initial data $(\eta_0, u_0, v_0)$ and the parameters $\delta$, $\varepsilon$, ${\rm R}$, and ${\rm W}$ satisfy Assumption \ref{assumption}, then we have the following estimates. 

\medskip

\noindent
I. Burgers equation

If the parameters $\delta$ and ${\rm W}$ and the initial data $\eta_0$ and $u_0$ satisfy
\begin{equation}\label{assumption I}
{\rm W_1}\le{\rm W}\le\delta^{-1}{\rm W}_2, 
 \quad |\eta_0|_{m+7}+\delta^{-1}\|(1+|D|_x)^{m+1}u_{0yy}\|\le M<\infty,
\end{equation}
then the following  error estimate holds.
\begin{equation}\label{error estimate I}
\mathcal{D}(t; \zeta^I, u^I, v^I, p^I)\le C\delta^2{\rm e}^{-c\varepsilon t}.
\end{equation}
II. Burgers equation with a fourth order dissipation term 

If the parameters $\delta$ and ${\rm W}$ and the initial data  $\eta_0$ and $u_0$ satisfy
\begin{equation}\label{assumption II}
{\rm W}=\delta^{-2}{\rm W}_2,\quad |\eta_0|_{m+12}+\delta^{-1}\|(1+|D|_x)^{m+1}u_{0yy}\|\le M<\infty,
\end{equation}
then the following  error estimate holds.
\begin{equation}\label{error estimate II}
 \mathcal{D}(t; \zeta^{II}, u^{II}, v^{II}, p^{II})\le C\delta^2{\rm e}^{-c\varepsilon t}.
\end{equation}
III. Burgers equation with  dispersion and nonlinear terms

If the parameters $\delta$ and ${\rm W}$ and the initial data $\eta_0$ and $u_0$ satisfy
\begin{equation}\label{assumption III}
{\rm W}_1\le{\rm W}\le{\rm W}_2, 
 \quad |\eta_0|_{m+13}+\delta^{-2}\|(1+|D|_x)^{m+1}(u_{0yy}-u_{yy}^{III}|_{t=0})\|\le M<\infty,
\end{equation}
then the following error estimate holds.
\begin{equation}\label{error estimate III}
 \mathcal{D}(t; \zeta^{III}, u^{III}, v^{III}, p^{III})\le C\delta^4{\rm e}^{-c\varepsilon t}.
\end{equation}
IV. Burgers equation with a fourth order dissipation, dispersion, and nonlinear terms

If the parameters $\delta$ and ${\rm W}$ and the initial data $\eta_0$ and $u_0$ satisfy
\begin{equation}\label{assumption IV}
{\rm W}=\delta^{-1}{\rm W}_2, 
 \quad |\eta_0|_{m+17}+\delta^{-2}\|(1+|D|_x)^{m+1}(u_{0yy}-u_{yy}^{IV}|_{t=0})\|\le M<\infty,
\end{equation}
then the following error estimate holds.
\begin{equation}\label{error estimate IV}
\mathcal{D}(t; \zeta^{IV}, u^{IV}, v^{IV}, p^{IV})\le C\delta^4{\rm e}^{-c\varepsilon t}.
\end{equation}
Here,  positive constants $C$ and $c$ depend on ${\rm R}_1, {\rm W}_1, {\rm W}_2, \alpha$, and $M$ but are independent of $\delta$, $\varepsilon$, ${\rm R}$, and ${\rm W}$. 
\end{thm}

\noindent
{\bf{Remark 2.2.}}
The assumptions for $u_{0yy}$ in \eqref{assumption I} and \eqref{assumption II} represent the restriction on the initial profile of the velocity. Moreover, the assumptions for $u_{0yy}$ in \eqref{assumption III} and \eqref{assumption IV} mean that the initial profile of the velocity have to be equal  to that of the approximate solution up to $O(\delta^2)$.

\medskip

\noindent
{\bf{Remark 2.3.}}
We see formally that the  order of error terms in \eqref{zeta I}  is of $O(\delta)$, which implies that the error estimates \eqref{error estimate I} and \eqref{error estimate II} are natural. In a similar way, we see that the error estimates \eqref{error estimate III} and \eqref{error estimate IV} are natural.

\medskip

\noindent
{\bf{Remark 2.4.}}
By introducing the slow time scale $\tau=\varepsilon t$, the norm decays exponentially and uniformly in $\tau$.

\medskip


\section{Approximate solutions}

In this section, following  Benney's perturbation method \cite{Benney} we will give approximate equations  by constructing approximate solutions. Hereafter, we assume $\varepsilon=\delta$. By straightforward calculation, we can rewrite \eqref{NS-2}--\eqref{BC0-2} as follows.

\begin{equation}\label{NS-3}
\left\{
 \begin{array}{lcl}
 \delta(u_t+\bar{u}u_x+\bar{u}_yv)+\dfrac2{\rm R}\delta p_x-\dfrac1{\rm R}(\delta^2 u_{xx}+u_{yy})=-\delta\dfrac2{\rm R}\eta u_{yy}+\delta^2f_1^{(2)}+\delta^3 f_1^{(3)}& \mbox{in} & \Omega, \; t>0, \\[3mm]
\delta^2(v_t+\bar{u}v_x)+\dfrac2{\rm R} p_y-\dfrac1{\rm R}\delta(\delta^2 v_{xx}+v_{yy})=\delta\dfrac2{\rm R}\eta p_{y}+\delta^2f_2^{(2)}+\delta^3 f_2^{(3)}& \mbox{in} & \Omega, \; t>0, \\[3mm]
  u_x+v_y=0 & \mbox{in} & \Omega, \; t>0,
 \end{array}
\right.
\end{equation}
\begin{equation}\label{BC1-3}
\left\{
 \begin{array}{lcl}
  \delta^2v_x+u_y-2(1+\delta\eta)^2
\eta=\delta^3h_1 & \mbox{on} & \Gamma, \; t>0, \\[1mm]
  p-\delta v_y-\dfrac{1}{\tan\alpha}\eta+\dfrac{\delta^2{\rm W}}{\sin\alpha}\eta_{xx}=\delta^2 h_2^{(2)}+\delta^3h^{(3)}_2  & \mbox{on} & \Gamma, \; t>0, 
\end{array}
\right.
\end{equation}
\begin{equation}\label{BC0-3}
u=v=0 \quad\mbox{on}\quad \Sigma, \; t>0,
\end{equation}
\begin{equation}\label{kinematic}
\eta_t+\eta_x-v=\delta^2h_3 \quad \mbox{on} \quad \Gamma, \; t>0,
\end{equation}
where 
\begin{equation}\label{definition of non-homogeneous term}
\begin{cases}
f_1^{(2)}=\dfrac{1}{\rm R}\big(3\eta^2u_{yy}-2\eta p_x+2y\eta_xp_y\big)+\eta_tu+y\eta_tu_y\\
\phantom{g^{(1)}=\dfrac{1}{\rm R}\big(}+y^2\eta_xu+2y(y-1)\eta u_x-y^2(y-2)\eta_x u_y
-u u_x-v u_y+2(2y-1)\eta v,\\[3mm]
f_2^{(2)}=\dfrac{1}{\rm R}\big(-2\eta^2 p_y+2\eta_xu_y+2\eta u_{xy}\big),\\[3mm]
h^{(2)}_2=2\eta\eta_x+\eta_x u+\eta u_x.
\end{cases}
\end{equation}

We proceed to derive the approximate equations following Benney \cite{Benney}.  Let $\eta=\eta(x, t)$ be a given function. For any $\delta\in(0, 1]$, let $(u, v, p)$ be the solution of \eqref{NS-3}--\eqref{BC0-3} and we expand $(u, v, p)$   as
\begin{equation}\label{perturbation expansion}
\begin{cases}
u=u_0+\delta u_1+\delta^2u_2+\cdots,\\
v=v_0+\delta v_1+\delta^2v_2+\cdots,\\
p=p_0+\delta p_1+\delta^2p_0+\cdots
\end{cases}
\end{equation}
and substitute these into \eqref{NS-3}--\eqref{BC0-3}, we obtain a sequence of perturbation equations for each order of $\delta$. Here, $u_0$ and $v_0$ are different from initial data defined in \eqref{initial} and hereafter we use this notation whenever it does not lead to confusion.
By assuming ${\rm W}=O(1)$, the $O(1)$, $O(\delta)$, and $O(\delta^2)$ problems are as follows.
\begin{equation}\label{order 1}
\begin{cases}
u_{0yy}=0,\quad p_{0y}=0,\quad u_{0x}+v_{0y}=0&\qquad {\rm in}\quad \Omega,\\
u_{0y}=2\eta,\quad p_0=\dfrac1{\tan\alpha}\eta&\qquad{\rm on}\quad \Gamma,\\
u_0=v_0=0&\qquad{\rm on}\quad \Sigma,
\end{cases}
\end{equation}
\begin{equation}\label{order delta1}
\begin{cases}
u_{1yy}={\rm R}(u_{0t}+(2y-y^2)u_{0x}+2(1-y)v_0)+2p_{0x}+2\eta u_{0yy}&\qquad {\rm in}\quad \Omega,\\
2p_{1y}=v_{0yy}+2\eta p_{0y}, \quad u_{1x}+v_{1y}=0&\qquad {\rm in}\quad \Omega,\\
u_{1y}=4\eta^2,\quad p_1=-u_{0x}&\qquad {\rm on}\quad \Gamma,\\
u_1=v_1=0&\qquad {\rm on}\quad \Sigma,
\end{cases}
\end{equation}
\begin{equation}\label{order delta2}
\begin{cases}
u_{2yy}={\rm R}(u_{1t}+(2y-y^2)u_{1x}+2(1-y)v_1)+2p_{1x}+2\eta u_{1yy}-u_{0xx}-{\rm R}f_1^{(2)}(\eta, u_0, v_0, p_0)&\quad {\rm in}\quad \Omega,\\
2p_{2y}=v_{1yy}+2\eta p_{1y}-{\rm R}\big(v_{0t}+(2y-y^2)v_{0x}\big)+{\rm R}f_2^{(2)}(\eta, u_0, v_0, p_0), \quad  u_{2x}+v_{2y}=0&\quad {\rm in}\quad \Omega,\\
u_{2y}=-v_{0x}+2\eta^3, \quad p_2=-u_{1x}+h^{(2)}_2(\eta, u_0)-\dfrac{{\rm W}}{\sin\alpha}\eta_{xx}&\quad {\rm on}\quad \Gamma,\\
u_2=v_2=0&\quad {\rm on}\quad \Sigma.
\end{cases}
\end{equation}
 Solving the above boundary value problem for the  ordinary differential equations, we have
\begin{equation}\label{definition of u0}
\begin{cases}
u_0=2y\eta,\\
v_0=-y^2\eta_x,\\
p_0=\dfrac1{\tan\alpha}\eta,
\end{cases}
\end{equation}
\begin{equation}\label{definition of u1}
\begin{cases}
u_1= \left(\frac{1}{3}y^3 -y\right){\rm R}\eta _t +\left\{ (y^2-2y)\frac{1}{\tan\alpha} +\left(\frac{1}{6}y^4-\frac{2}{3}y\right){\rm R}\right\}\eta _x+4y\eta^2 ,\\
v_1= \left(-\frac{1}{12}y^4+\frac{1}{2}y^2\right){\rm R}\eta _{xt}+\left\{ \big(-\frac{1}{3}y^3+y^2\big)\frac{1}{\tan\alpha} +\left(-\frac{1}{30}y^5+\frac{1}{3}y^2\right){\rm R}\right\}\eta _{xx}-4y^2\eta\eta_x, \\
p_1= -(1+y)\eta _x,
\end{cases}
\end{equation}
\begin{equation}\label{definition of u2}
\begin{cases}
u_2=\left( \frac{1}{60}y^5-\frac{1}{6}y^3+\frac{5}{12}y\right) {\rm R}^2\eta _{tt}\\
 \qquad\quad+  \big\{ \left( \frac{1}{12}y^4-\frac{1}{3}y^3+\frac{2}{3}y\right) \frac{{\rm R}}{\tan\alpha}+\left( -\frac{1}{252}y^7+\frac{1}{45}y^6-\frac{1}{12}y^4-\frac{1}{9}y^3+\frac{101}{180}y\right) {\rm R}^2 \big\}\eta _{xt}\\
 \qquad\quad+ \big\{ \left( -\frac{2}{3}y^3-y^2+5y\right) +  \left( -\frac{1}{90}y^6+\frac{1}{15}y^5-\frac{1}{6}y^4+\frac{2}{5}y\right) \frac{{\rm R}}{\tan\alpha} \\
 \qquad\qquad  +\left( -\frac{1}{560}y^8+\frac{2}{315}y^7-\frac{1}{18}y^4+\frac{121}{630}y\right) {\rm R}^2 \big\}\eta _{xx}\\
\qquad\quad+2y\eta^3+{\rm R}(\frac43y^3-4y)\eta\eta_t+\{{\rm R}(y^4-4y)+(3y^2-6y)\frac1{\tan\alpha}\}\eta\eta_x, \\
v_2=\left( -\frac{1}{360}y^6+\frac{1}{24}y^4-\frac{5}{24}y^2\right) {\rm R}^2\eta _{xtt}\\
 \qquad\quad +\left\{ \left( -\frac{1}{60}y^5+\frac{1}{12}y^4-\frac{1}{3}y^2\right) \frac{{\rm R}}{\tan\alpha} 
+\left( \frac{1}{2016}y^8-\frac{1}{315}y^7+\frac{1}{60}y^5+\frac{1}{36}y^4-\frac{101}{360}y^2\right) {\rm R}^2\right\} \eta _{xxt}\\
 \qquad\quad+\big \{ \left( \frac{1}{6}y^4+\frac{1}{3}y^3-\frac{5}{2}y^2\right) +\left( \frac{1}{630}y^7-\frac{1}{90}y^6+\frac{1}{30}y^5-\frac{1}{5}y^2\right) \frac{{\rm R}}{\tan\alpha} \\
\qquad\qquad+\left( \frac{1}{5040}y^9-\frac{1}{1260}y^8+\frac{1}{90}y^5-\frac{121}{1260}y^2\right) {\rm R}^2\big \}\eta _{xxx}\\
\qquad\quad-3y^2\eta^2\eta_x+{\rm R}\big(-\frac13y^4+2y^2\big)(\eta_x\eta_t+\eta\eta_{tx})+\big\{{\rm R}\big(-\frac15y^5+2y^2\big)+(-y^3+3y^2)\frac1{\tan\alpha}\big\}(\eta_x^2+\eta\eta_{xx}),\\
p_2=\left( \frac{1}{2}y+\frac{1}{6}\right) {\rm R}\eta _{xt}
+\left \{-\frac{\rm W}{\sin\alpha} +\left(-\frac{1}{2}y^2+y+\frac{1}{2}\right)\frac{1}{\tan\alpha}  +\left( -\frac{1}{10}y^5+\frac{1}{6}y^4+\frac{1}{3}y+\frac{1}{10}\right) {\rm R}\right \} \eta _{xx}\\
\qquad\quad +\big\{{\rm R}(4y-4)-5y+3\big\}\eta\eta_x.
\end{cases}
\end{equation}
Using the above expressions, we put
\begin{equation}\label{pre approximate solution III}
\begin{cases}
u^{III}_{0}(y;\eta):=u_0, \quad v^{III}_0(y;\eta):=v_0, \quad p^{III}_0(y;\eta):=p_0,\\
u^{III}_1(y;\eta):=u_1, \quad v^{III}_1(y;\eta):=v_1, \quad p^{III}_1(y;\eta):=p_1,\\
u^{III}_2(y;\eta):=u_2, \quad v^{III}_2(y;\eta):=v_2, \quad p^{III}_2(y;\eta):=p_2.
\end{cases}
\end{equation}
In view of the perturbation expansion \eqref{perturbation expansion} substituting $v=v_0^{III}+\delta v_1^{III}+\delta^2 v_2^{III}$  into \eqref{kinematic}, we obtain the approximate equation 
\[
\eta_t+2\eta_x+4\varepsilon\eta\eta_x-\frac{8}{15}({\rm R}_c-{\rm R})\delta\eta_{xx}+C_1\delta^2\eta_{xxx}+C_2\varepsilon\delta(\eta\eta_{xx}+\eta_x^2)
+2\varepsilon^2\eta^2\eta_x=O(\delta^3),
\]
where $C_1=2+\frac{32}{63}{\rm R}^2-\frac{40}{63}\frac{\rm R}{\tan\alpha}$ and $C_2=\frac{16}{5}{\rm R}-\frac2{\tan\alpha}$.

Thus far we have assumed ${\rm W}=O(1)$. 
Taking into account that $\rm W$ is contained only in the second equation in \eqref{BC1-3} and modifying  $O(\delta)$ problem under  the assumption ${\rm W}\le O(\delta^{-1})$, we see that $(u_0^{I}, v_0^{I}, p_0^{I})$ and $(u_1^{I}, v_1^{I}, p_1^{I})$, which are defined by
\begin{equation}\label{pre approximate solution I}
\begin{cases}
u^{I}_{0}(y;\eta):=u_0, \quad v^{I}_0(y;\eta):=v_0, \quad p^{I}_0(y;\eta):=p_0,\\
u^{I}_1(y;\eta):=u_1, \quad v^{I}_1(y;\eta):=v_1, \quad p^{I}_1(y;\eta):=p_1-\frac{\delta{\rm W}}{\sin\alpha}\eta_{xx},
\end{cases}
\end{equation}
are the solutions of the problems. Putting $v=v_0^{I}+\delta v_1^{I}$ and substituting this into \eqref{kinematic}, we obtain the approximate equation
\[
\eta_t+2\eta_x+4\varepsilon\eta\eta_x-\frac{8}{15}({\rm R}_c-{\rm R})\delta\eta_{xx}=O(\delta^2).
\]
Similarly, modifying  $O(1)$ and $O(\delta)$ problems under  the assumption ${\rm W}= O(\delta^{-2})$ and putting
\begin{equation}\label{pre approximate solution II}
\begin{cases}
u^{II}_{0}:=u_0, \quad v^{II}_0:=v_0, \quad p^{II}_0:=p_0-\frac{\delta^2{\rm W}}{\sin\alpha}\eta_{xx},\\
u^{II}_1:=u_1-\frac{\delta^2{\rm W}}{\sin\alpha}(y^2-2y)\eta_{xxx}, \quad v^{II}_1:=v_1+\frac{\delta^2{\rm W}}{\sin\alpha}\big(\frac13y^3-y^2\big)\eta_{xxxx}, \quad p^{II}_1:=p_1,
\end{cases}
\end{equation}
we obtain the approximate equation
\[
\eta_t+2\eta_x+4\varepsilon\eta\eta_x-\frac{8}{15}({\rm R}_c-{\rm R})\delta\eta_{xx}+\frac23\frac{{\rm W}_2}{\sin\alpha}\delta\eta_{xxxx}=O(\delta^2).
\]
Moreover, putting 
\begin{equation}\label{pre approximate solution IV}
\begin{cases}
u^{IV}_{0}:=u_0, \quad v^{IV}_0:=v_0, \quad p^{IV}_0:=p_0,\\
u^{IV}_1:=u_1, \quad v^{IV}_1:=v_1, \quad p^{IV}_1:=p_1-\frac{\delta{\rm W}}{\sin\alpha}\eta_{xx},\\
 u^{IV}_2:=u_2-\frac{\delta{\rm W}}{\sin\alpha}(y^2-2y)\eta_{xxx}, \quad v^{IV}_2:=v_2+\frac{\delta{\rm W}}{\sin\alpha}\big(\frac13y^3-y^2\big)\eta_{xxxx}, \quad p^{IV}_2:=p_2+\frac{{\rm W}}{\sin\alpha}\eta_{xx}
\end{cases}
\end{equation}
and $v=v_0^{IV}+\delta v_1^{IV}+\delta^2 v_2^{IV}$ and substituting this into \eqref{kinematic}, we obtain the approximate equation
\begin{align*}
&\eta_t+2\eta_x+4\varepsilon\eta\eta_x-\frac{8}{15}({\rm R}_c-{\rm R})\delta\eta_{xx}\\
&\phantom{\eta_t}+C_1\delta^2\eta_{xxx}+C_2\varepsilon\delta(\eta\eta_{xx}+\eta_x^2)+2\varepsilon^2\eta^2\eta_x+\frac23\frac{{\rm W}_2}{\sin\alpha}\delta^2\eta_{xxxx}=O(\delta^3)
\end{align*}
under the assumption ${\rm W}= O(\delta^{-1})$.


\section{Energy estimate}

In this section, we will derive energy estimates. Let $\eta=\eta(x,t)$ be a fixed function and $(u, v, p)=(u(y;\eta), v(y;\eta), p(y;\eta))$ be an approximate solution constructed from $\eta$ satisfying $u_x+v_y=0$ and $u|_{y=0}=v|_{y=0}=0$, which will be defined in the next section. 
Using the approximate solution, we define $\psi_1, \psi_2, \phi_1, \phi_2, \phi_3$ by the following equality.
\begin{equation}\label{definition of phi}
\begin{cases}
\psi_1(y;\eta):=\dfrac1{\delta^3}\bigg\{\delta(u_t+\bar{u}u_x+\bar{u}_yv)+\dfrac2{\rm R}\delta p_x-\dfrac1{\rm R}(\delta^2 u_{xx}+u_{yy})-\delta f_1^{(1)}(\eta, u, v, p)\bigg\}, \\[3mm]
\psi_2(y;\eta):=\dfrac1{\delta^3}\bigg\{\delta^2(v_t+\bar{u}v_x)+\dfrac2{\rm R} p_y-\dfrac1{\rm R}\delta(\delta^2 v_{xx}+v_{yy})-\delta f_2^{(1)}(\eta, u,  p)\bigg\}, \\[3mm]
\phi_1(\eta):=\dfrac1{\delta^3}\big\{\delta^2v_x+u_y-2(1+\delta\eta)^2\eta\big\}|_{y=1},\\
\phi_2(\eta):=\dfrac1{\delta^3}\bigg\{p-\delta v_y-\dfrac{1}{\tan\alpha}\eta+\dfrac{\delta^2{\rm W}}{\sin\alpha}\eta_{xx}-\delta^2h^{(2)}_2(\eta, u)\bigg\}\bigg|_{y=1},\\
\phi_3(\eta):=\dfrac1{\delta^3}\{\eta_t+\eta_x-v-\delta^2h_3(\eta)\}|_{y=1},
\end{cases}
\end{equation}
where
\begin{equation}\label{tildeg}
f_1^{(1)}=-\dfrac2{\rm R}\eta u_{yy}+\delta f_1^{(2)},\quad f_2^{(1)}=\dfrac2{\rm R}\eta p_{y}+\delta f_2^{(2)}.
\end{equation}
Then the approximate solution satisfies the following equations.
\begin{equation}\label{NSapp}
\left\{
 \begin{array}{lcl}
 \delta(u_t+\bar{u}u_x+\bar{u}_yv)+\dfrac2{\rm R}\delta p_x-\dfrac1{\rm R}(\delta^2 u_{xx}+u_{yy})=\delta f_1^{(1)}(\eta, u, v, p)+\delta^3\psi_1(y;\eta)& \mbox{in} & \Omega, \; t>0, \\[3mm]
\delta^2(v_t+\bar{u}v_x)+\dfrac2{\rm R} p_y-\dfrac1{\rm R}\delta(\delta^2 v_{xx}+v_{yy})=\delta f_2^{(1)}(\eta, u, p)+\delta^3 \psi_2(y;\eta)& \mbox{in} & \Omega, \; t>0, \\[3mm]
  u_x+v_y=0 & \mbox{in} & \Omega, \; t>0,
 \end{array}
\right.
\end{equation}
\begin{equation}\label{BC1app}
\left\{
 \begin{array}{lcl}
  \delta^2v_x+u_y-2(1+\delta\eta)^2\eta=\delta^3\phi_1(\eta) & \mbox{on} & \Gamma, \; t>0, \\[1mm]
  p-\delta v_y-\dfrac{1}{\tan\alpha}\eta+\dfrac{\delta^2{\rm W}}{\sin\alpha}\eta_{xx}=\delta^2h_2^{(2)}(\eta, u)+\delta^3\phi_2(\eta)
   & \mbox{on} & \Gamma, \; t>0, \\[2mm]
  \eta_t+\eta_x-v=\delta^2h_3(\eta)+\delta^3\phi_3(\eta)& \mbox{on} & \Gamma, \; t>0,
 \end{array}
\right.
\end{equation}
\begin{equation}\label{BC0app}
u=v=0 \quad\mbox{on}\quad \Sigma, \; t>0.
\end{equation}
Note that in the next section, we will give  explicit forms of $\psi_1$, $\psi_2$, $\phi_1$, $\phi_2$, and $\phi_3$. Let $(\eta^\delta, u^\delta, v^\delta, p^\delta)$ be the solution of \eqref{NS-3}--\eqref{BC0-3} and we set
\[
H:=\eta^\delta-\eta,\quad U:=u^\delta-u,\quad V:=v^\delta-v, \quad P:=p^\delta-p.
\]
Taking the difference between  \eqref{NS-3}--\eqref{kinematic} and \eqref{NSapp}--\eqref{BC0app}, we have
\begin{equation}\label{NSerror}
\left\{
\begin{array}{lcl}
 \delta(U_t+\bar{u}U_x+\bar{u}_yV)+\dfrac2{\rm R}\delta P_x-\dfrac1{\rm R}(\delta^2 U_{xx}+U_{yy})\\[3mm]
\hspace{4cm}=F_1+\delta^3 f_1^{(3)}(\eta^\delta, u^\delta, v^\delta, p^\delta)-\delta^3\psi_1(y;\eta)& \mbox{in} & \Omega, \; t>0, \\
\delta^2(V_t+\bar{u}V_x)+\dfrac2{\rm R} P_y-\dfrac1{\rm R}\delta(\delta^2 V_{xx}+V_{yy})\\[3mm]
\hspace{4cm}=F_2+\delta^3 f_2^{(3)}(\eta^\delta, u^\delta, v^\delta, p^\delta)-\delta^3\psi_2(y;\eta)& \mbox{in} & \Omega, \; t>0, \\
  U_x+V_y=0 & \mbox{in} & \Omega, \; t>0,
\end{array}
\right.
\end{equation}
\begin{equation}\label{BC1error}
\left\{
 \begin{array}{lcl}
  \delta^2V_x+U_y-\big(2+b(\eta^\delta, \eta)\big)H=\delta^3h_1(\eta^\delta, u^\delta, v^\delta)-\delta^3\phi_1(\eta) & \mbox{on} & \Gamma, \; t>0, \\[1mm]
  P-\delta V_y-\dfrac{1}{\tan\alpha}H+\dfrac{\delta^2{\rm W}}{\sin\alpha}H_{xx}=G_2
+\delta^3h_2^{(3)}(\eta^\delta, u^\delta, v^\delta)-\delta^3\phi_2(\eta)
   & \mbox{on} & \Gamma, \; t>0, \\[2mm]
  H_t+H_x-V= G_3-\delta^3\phi_3(\eta)
 & \mbox{on} & \Gamma, \; t>0,
 \end{array}
\right.
\end{equation}
\begin{equation}\label{BC0error}
U=V=0 \quad\mbox{on}\quad \Sigma, \; t>0,
\end{equation}
where 
\begin{equation}\label{definition of error}
\begin{cases}
F_1= \delta\big(f_1^{(1)}(\eta^\delta, u^\delta, v^\delta, p^\delta)-f_1^{(1)}(\eta, u, v, p)\big), \quad F_2= \delta\big(f_2^{(1)}(\eta^\delta, u^\delta, p^\delta)-f_2^{(1)}(\eta, u, p)\big),\\ b=2\delta\big(\delta(\eta^\delta)^2+(2+\delta\eta)\eta^\delta+\delta\eta^2+2\eta\big),\\
G_2=\delta^2\big(h^{(2)}_2(\eta^\delta, u^\delta, v^\delta)-h^{(2)}_2(\eta, u, v)\big), \quad G_3=\delta^2\big(h_3(\eta^\delta)-h_3(\eta)\big).
\end{cases}
\end{equation}
For convenience, we set
\[
\bm{U}:=(U, \delta V)^{\rm T},\quad \bm{F}:=(F_1, F_2)^{\rm T},\quad \bm{f}:=(f_1^{(3)}, f_2^{(3)})^{\rm T}, \quad \bm{\psi}:=(\psi_1, \psi_2)^{\rm T}.
\]

We proceed to derive an energy estimate to \eqref{NSerror}--\eqref{BC0error} following \cite{Ueno}.
In view of the energies obtained in \cite{Ueno} (see (3.6)--(3.8) and (3.24) in \cite{Ueno}), we put 
\begin{align*}
&\mathscr {E}_0(H, \bm{U})
:= \delta^2\| V\|^2+\dfrac{2}{\rm R}\biggl(\frac{1}{\tan\alpha}|H|_0^2
  +\frac{\delta^2{\rm W}}{\sin\alpha}|H_x|_0^2\biggr) \\
&\phantom{\mathscr {E}_0(H, \bm{U})
:=}
 +\beta_1\bigg\{\delta^2\|\bm{U}_x\|^2
  +\frac{2}{\rm R}\biggl(\dfrac{1}{\tan\alpha}\delta^2|H_x|_0^2
  +\dfrac{\delta^2{\rm W}}{\sin\alpha}\delta^2|H_{xx}|_0^2\biggr)\bigg\} \notag \\
&\phantom{\mathscr {E}_0(H, \bm{U})
:=}
 +\beta_2\bigg\{\delta^4\|\bm{U}_{xx}\|^2
  +\frac{2}{\rm R}\biggl(\dfrac{1}{\tan\alpha}\delta^4|H_{xx}|_0^2
  +\dfrac{\delta^2{\rm W}}{\sin\alpha}\delta^4|H_{xxx}|_0^2\biggr)\bigg\} \notag \\
&\phantom{\mathscr {E}_0(H, \bm{U})
:=}
 +\beta_3\bigg\{\delta^2\|\bm{U}_t\|^2
  +\frac{2}{\rm R}\biggl(\dfrac{1}{\tan\alpha}\delta^2|H_t|_0^2
  +\dfrac{\delta^2{\rm W}}{\sin\alpha}\delta^2|H_{tx}|_0^2\biggr)\bigg\},\notag\\
&\mathscr{F}_0(H, \bm{U}, P)
:= \delta\|\bm{U}_x\|^2+\delta\|\partial_y^{-1} P_x\|^2
 +\delta|H_x|_0^2
  +\delta^3{\rm W}|H_{xx}|_0^2
  +\delta^5{\rm W}^2|H_{xxx}|_0^2\notag \\
&\phantom{\mathscr{F}_0(H, \bm{U}, P)
:= }
 +\delta\|\nabla_\delta\bm{U}_x\|^2
  +\delta^3\|\nabla_\delta\bm{U}_{xx}\|^2
  +\delta\|\nabla_\delta\bm{U}_t\|^2. \notag
\end{align*}
Here, $\beta_1, \beta_2, $ and $\beta_3$ are appropriate positive constants (see (3.28) in \cite{Ueno}). 
Integrating by parts and using the third equation in \eqref{BC1error} and Poincar\'e's inequality, we see that for any $\epsilon>0$ there exists a positive constant $C_\epsilon$ such that
\begin{align*}
&\delta^3|(\{\bm{F}+\delta^3\bm{f}-\delta^3\bm{\psi}\}_{xx}, \bm{U}_{xx})_\Omega|\le\epsilon\delta^5\|\bm{U}_{xxx}\|^2+C_\epsilon\delta(\|\bm{F}_x\|^2+\delta^6\|\bm{f}_x\|^2+\delta^6\|\bm{\psi}_x\|^2),\\
&|(H, (bH)_x)_\Gamma|\le\epsilon\delta|H_x|_0^2+C_\epsilon\delta^{-1}|(bH)_x|_0^2,\\
&\delta^2{\rm W}|(H_{xx}, (bH)_x)_\Gamma|\le\epsilon\delta^3{\rm W}|H_{xx}|_0^2+C_\epsilon\delta{\rm W}|(bH)_x|_0^2,\\
&\delta^2{\rm W}|(H_{xx}, G_3-\delta^3\phi_3)_\Gamma|\le\epsilon\delta^3{\rm W}|H_{xx}|_0^2+C_\epsilon\delta{\rm W}(|G_{3}|_0^2+\delta^6|\phi_{3}|_0^2),\\
&\delta^6{\rm W}|(H_{xxxx}, \delta^3\phi_{3xx})_\Gamma|\le\epsilon\delta^5{\rm W}^2|H_{xxx}|_0^2+C_\epsilon\delta^{13}|\phi_{3xxx}|_0^2,\\
&\delta^4{\rm W}|(H_{xxt},  G_{3t}-\delta^3\phi_{3t})_\Gamma|
  \le\epsilon(\delta^5{\rm W}^2|H_{xxx}|_0^2+\delta^5\|U_{xxx}\|_0^2)\\
&\phantom{\delta^4{\rm W}|(H_{xxt}, }+C_\epsilon(1+{\rm W}^2)\delta^3(|G_{3t}|_0^2+\delta^6|\phi_{3t}|_0^2)+\delta^5(|G_{3xx}|_0^2+\delta^6|\phi_{3xx}|_0^2).
\end{align*}
Here, we used the inequality $|V(\cdot,1)|_0=|V(\cdot,1)-V(\cdot,0)|_0\leq\|V_y\|=\|U_x\|$ 
thanks to the third equation in \eqref{NSerror} and the second equation in \eqref{BC0error}. 
In the following, we use frequently this type of inequality without any comment. 
Taking into account the above inequality and (3.27) in \cite{Ueno}, we need to estimate the following quantities.
\begin{align}
&\mathscr{N}^{1}_0(Z_1)
:= (\delta{\rm W}+\delta^{-1})|(bH)_x|_0^2+\delta^3|(bH)_{xx}|_0^2+\delta|(bH)_t|_0^2\label{definition of N1}\\
&\phantom{\mathscr{N}^{1}_0(Z_1)
:=}
+\delta^{-1}|G_2|_0^2+\delta|G_{2x}|_0^2  +\delta^2||D_x|^{\frac12}G_{2x}|_0^2
 +\delta|(G_{2t}, \delta V_t)_\Gamma|  \notag\\
&\phantom{\mathscr{N}^{1}_0(Z_1)
:=}
+\delta{\rm W}|G_3|_0^2+\delta^3|G_{3x}|_0^2+\delta^5|G_{3xx}|_0^2+\delta^3{\rm W}^2|G_{3t}|_0^2+\delta^6{\rm W}|(H_{xxxx}, G_{3xx})_\Gamma|\notag\\
&\phantom{\mathscr{N}^{1}_0(Z_1)
:=}
 +\delta^{-1}\|\bm{F}\|^2 +\delta\|\bm{F}_{x}\|^2+\delta|(\bm{F}_t, \bm{U}_t)_\Omega|,\notag\\
&\mathscr{N}_0^{2}(Z_2)
:= \delta^{5}|h_1|_0^2+\delta^7|h_{1x}|_0^2+\delta^8||D_x|^{\frac12}h_{1x}|_0^2+\delta^4|(h_{1t}, U_t)_\Gamma|\label{definition of N2}\\
&\phantom{\mathscr{N}_0^{2}(Z_2)
:= }+\delta^{5}|h^{(3)}_2|_0^2+\delta^7|h^{(3)}_{2x}|_0^2+\delta^8||D_x|^{\frac12}h^{(3)}_{2x}|_0^2+\delta^4|(h^{(3)}_{2t}, \delta V_t)_\Gamma|\notag\\
&\phantom{\mathscr{N}_0^{2}(Z_2)
:= }+\delta^{5}\|\bm{f}\|^2+\delta^7\|\bm{f}_{x}\|^2+\delta^4|(\bm{f}_t, \bm{U}_t)_\Omega|,\notag\\
&\mathscr{N}^{3}_0(Z_3)
:= \delta^{5}|\phi_1|_0^2+\delta^7|\phi_{1x}|_0^2+\delta^8||D_x|^{\frac12}\phi_{1x}|_0^2+\delta^7|\phi_{1t}|_0^2+\delta^{5}|\phi_2|_0^2+\delta^7|\phi_{2x}|_0^2\label{definition of N3}\\
&\phantom{\mathscr{N}^{3}_0(Z_3)
:= }+\delta^8||D_x|^{\frac12}\phi_{2x}|_0^2
+\delta^7|\phi_{2t}|_0^2 
 +\delta^{7}{\rm W}|\phi_3|_0^2+\delta^9|\phi_{3x}|_0^2+\delta^{11}|\phi_{3xx}|_0^2\notag\\
&\phantom{\mathscr{N}^{3}_0(Z_3)
:= }+\delta^{13}|\phi_{3xxx}|_0^2 +\delta^9{\rm W}^2|\phi_{3t}|_0^2\notag 
 +\delta^{5}\|\bm{\psi}\|^2+\delta^7\|\bm{\psi}_{x}\|^2+\delta^7\|\bm{\psi}_t\|^2,\notag
\end{align}
where
\[
Z_1=(H, \bm{U}, bH, G_2, G_3, \bm{F}), \quad Z_2=(\bm{U}, h_1, h^{(3)}_2, h_3, \bm{f}),\quad Z_3=(\phi_1, \phi_2, \phi_3, \bm{\psi}).
\]
For an integer $m\ge2$, we set 
\begin{align}
&\mathscr{E}_m(H, \bm{U}) := \sum_{k=0}^m \mathscr{E}_0(\partial_x^kH, \partial_x^k\bm{U}), \quad 
\mathscr{F}_m(H, \bm{U}, P) := \sum_{k=0}^m \mathscr{F}_0(\partial_x^kH, \partial_x^k\bm{U}, \partial_x^k P),\label{Em}\\
&\mathscr{N}_m^{1}(H, \bm{U}, P;\eta):= \sum_{k=0}^m \big\{\mathscr{N}^{1}_0(\partial_x^k Z_1)+|(\partial_x^kH,\partial_x^k G_3)_\Gamma|\big\},\label{definition of Nm1}\\
&\mathscr{N}_m^2(\bm{U}):= \sum_{k=0}^m \mathscr{N}^{2}_0(\partial_x^k Z_2),\label{definition of Nm2}\\
&\mathscr{N}_m^{3}(H;\eta):= \sum_{k=0}^m \big\{\mathscr{N}^{3}_0(\partial_x^k Z_3)+|(\partial_x^kH,\delta^3\partial_x^k \phi_3)_\Gamma|\big\}.\label{definition of Nm3}
\end{align}
Here, the terms $\sum_{k=0}^m|(\partial_x^k H, \partial_x^kG_3)_\Gamma|$ and $\sum_{k=0}^m|(\partial_x^kH,\delta^3\partial_x^k \phi_3)_\Gamma|$ come from (3.30) in \cite{Ueno}.
Applying $\partial_x^k$ to \eqref{NSerror}--\eqref{BC0error}, using [16, Proposition 3.2], 
and adding the resulting inequalities for $0\le k\le m$, we obtain the following lemma.

\begin{lem}\label{lem energy estimate NS}
There exist small positive constants ${\rm R}_0$ and  $\alpha_0$ such that if $0<{\rm R_1}\le{\rm R}\le{\rm R}_0$, ${\rm W}_1\le{\rm W}$, and $0<\alpha\le\alpha_0$, 
then the solution $(H, U, V, P)$ of \eqref{NSerror}--\eqref{BC0error} satisfies 
\begin{equation}\label{energy estimate NS}
\frac{{\rm d}}{{\rm d}t}\mathscr{E}_m+\mathscr{F}_m\le C( \mathscr{N}_m^1+\mathscr{N}_m^2+\mathscr{N}_m^3),
\end{equation}
where the constant $C$ is independent of  $\delta$, ${\rm R}$,  and ${\rm W}$. 
\end{lem}

\medskip
For later use, we modify the energy and the dissipation functions $\mathscr{E}_m$ and $\mathscr{F}_m$ as 
\begin{align}
&\tilde{\mathscr{E}}_m(H, \bm{U}) := \mathscr{E}_m(H, \bm{U})+\|(1+|D_x|)^m U\|^2+\|(1+|D_x|)^m U_y\|^2, \label{tildeEm}\\
&\tilde{\mathscr{F}}_m(H, \bm{U}, P) := \mathscr{F}_m(H, \bm{U}, P)+
 \delta|(1+\delta|D_x|)^{\frac52}H_t|_m^2+(\delta^2{\rm W})^2\delta^2||D_x|^{\frac72}H|_m^2 \label{Fm}\\
&\phantom{\tilde{\mathscr{F}}_m(H, \bm{U}, P) := }
 +\delta^{-1}\|(1+|D_x|)^m(1+\delta|D_x|)(\nabla_\delta P, U_{yy})\|^2
 +\delta\|(1+|D_x|)^{m-1}\nabla_\delta P_t\|^2. \notag
\end{align}
We also introduce another energy function $\mathscr{D}_m$ by
\begin{align}
&\mathscr{D}_m(H, \bm{U}) := |(1+\delta|D_x|)^2H|_m^2+\delta^2\|(1+|D_x|)^mV\|^2+\delta^2\|(1+|D_x|)^m\bm{U}_x\|^2 \label{definition of D}\\
&\phantom{\mathscr{D}_m(H, \bm{U}) =}+\|(1+|D_x|)^mD_{\delta}^2\bm{U}\|^2+(\delta^2{\rm W})^2|(1+\delta|D_x|)H_x|_{m+1}^2+\sqrt{\delta^2{\rm W}}\|(1+|D_x|)^m\delta V_{xy}\|^2,\notag
\end{align}
which does not include any time derivatives. Setting $\tilde{E}_m=\tilde{\mathscr{E}}_m(\eta^\delta, \bm{u}^\delta)$ and  $\tilde{F}_m=\tilde{\mathscr{F}}_m(\eta^\delta, \bm{u}^\delta, p^\delta)$ 
and using [16, Theorem 2.2 and Proposition 6.1],  the following uniform estimate holds.
\begin{prop}\label{uniform estimate}
There exist small positive constants ${\rm R}_0$ and $\alpha_0$ such that the following statement holds: 
Let $m$ be an integer satisfying $m\ge2$, $0<{\rm R}_1\le {\rm R}_0$, $0<{\rm W}_1\le{\rm W}_2$, and $0<\alpha\le \alpha_0$. There exists small positive constant $c_0$ such that
if the initial data $(\eta_0, u_0, v_0)$ and the parameters $\delta$, $\varepsilon$, ${\rm R}$, and ${\rm W}$ satisfy Assumption \ref{assumption} and ${\rm W}\le\delta^{-2}{\rm W}_2$, then the solution $(\eta^\delta, u^\delta, v^\delta, p^\delta)$ of \eqref{NS-2}--\eqref{initial} satisfies
\[
\tilde{E}_2(t)\le c_0, \quad \sup_{t\ge0}\tilde{E}_{m+1}(t)+\int_0^\infty\tilde{F}_{m+1}(t){\rm d}t\le C,\quad \tilde{E}_{m+1}(t)\le C{\rm e}^{-c\delta t}.
\]
Here,  positive constants $C$ and $c$ depend on ${\rm R}_1, {\rm W}_1, {\rm W}_2, \alpha$, and $M$ but are independent of $\delta$, $\varepsilon$, ${\rm R}$, and ${\rm W}$. 
\end{prop}
\bigskip

Moreover, we easily obtain the following lemma.
\begin{lem}\label{energy estimate app}
Let $\alpha>0$, $0<{\rm R}_1\le{\rm R}<{\rm R}_c$. There exists small positive constant  $c_1$  such that if $s\ge2$ and $|\eta_0|_2^2\le c_1$, then the problems \eqref{zeta I}--\eqref{zeta IV} under the initial condition $\zeta|_{\tau=0}=\eta_0$ have  unique solutions $\zeta^I$, $\zeta^{II}$, $\zeta^{III}$, and $\zeta^{IV}$, respectively, which satisfy 
\begin{align*}
&\sup_{\tau\ge0}|\zeta^I(\tau)|_s^2+\int_0^\infty |\zeta^I_x(\tau)|_s^2{\rm d}\tau\le C|\eta_0|_s^2,\quad |\zeta^I(\tau)|_s^2\le C|\eta_0|_s^2{\rm e}^{-c\delta t},\\
&\sup_{\tau\ge0}|\zeta^{II}(\tau)|_s^2+\int_0^\infty \big(|\zeta^{II}_x(\tau)|_s^2+|\zeta^{II}_{xx}(\tau)|_s^2\big){\rm d}\tau\le C|\eta_0|_s^2 ,\quad |\zeta^{II}(\tau)|_s^2\le C|\eta_0|_s^2{\rm e}^{-c\delta t}, \\
&\sup_{\tau\ge0}|\zeta^{III}(\tau)|_s^2+\int_0^\infty |\zeta^{III}_x(\tau)|_s^2{\rm d}\tau\le C|\eta_0|_s^2,\quad |\zeta^{III}(\tau)|_s^2\le C|\eta_0|_s^2{\rm e}^{-c\delta t},\\
&\sup_{\tau\ge0}|\zeta^{IV}(\tau)|_s^2+\int_0^\infty \big(|\zeta^{IV}_x(\tau)|_s^2+\delta|\zeta^{IV}_{xx}(\tau)|_s^2\big){\rm d}\tau\le C|\eta_0|_s^2,\quad |\zeta^{IV}(\tau)|_s^2\le C|\eta_0|_s^2{\rm e}^{-c\delta t}. 
\end{align*}
Here, ${\rm R}_c=\frac54\frac1{\tan\alpha}$ is the critical Reynolds number and positive constants $C$ and $c$ are independent of $\delta$ and ${\rm R}$.
\end{lem}


\section{Error estimate}

We will show \eqref{error estimate III} under  Assumption \ref{assumption} and  \eqref{assumption III}.
We can show the other claims in Theorem \ref{main theorem} in the same way as the proof of \eqref{error estimate III} and we will comment about the difference at the end of this section.
Let $\zeta^{III}$ be the solution of \eqref{zeta III} under the initial condition $\zeta^{III}|_{\tau=0}=\eta_0$ and we put $\eta^{III}(x, t):=\zeta^{III}(x-2t, \varepsilon t)$ and 
\begin{equation}\label{approximate solution III}
\begin{cases}
u^{III}(x,y,t):=u_0^{III}(y;\eta^{III}(x,t))+\delta u^{III}_1(y;\eta^{III}(x,t))+\delta^2 u^{III}_2(y;\eta^{III}(x,t)),\\
v^{III}(x,y,t):=u^{III}_0(y;\eta^{III}(x,t))+\delta v^{III}_1(y;\eta^{III}(x,t))+\delta^2 v^{III}_2(y;\eta^{III}(x,t)),\\
p^{III}(x,y,t):=p^{III}_0(y;\eta^{III}(x,t))+\delta p^{III}_1(y;\eta^{III}(x,t))+\delta^2 p^{III}_2(y;\eta^{III}(x,t)),
\end{cases}
\end{equation}
where $u^{III}_0, v^{III}_0, p^{III}_0, \ldots$ were defined by \eqref{definition of u0}--\eqref{pre approximate solution III}.
Then, we have
\begin{align}\label{t-derivative of eta III}
\eta_t^{III}=&-2\eta^{III}_x+\dfrac{8}{15}({\rm{R}}_c-{\rm{R}})\delta\eta^{III}_{xx}-C_1\delta^2\eta^{III}_{xxx}\\
&-4\delta\eta^{III}\eta^{III}_x -\delta^2\big\{C_2\big(\eta^{III}\eta^{III}_{xx}+(\eta^{III}_x)^2\big)+2(\eta^{III})^2\eta^{III}_x\big\}.\notag
\end{align}
Using the approximate solutions \eqref{approximate solution III}, we define $\psi_1$, $\psi_2$, $\phi_1$, $\phi_2$, and $\phi_3$ by \eqref{definition of phi}.
By using the equality \eqref{t-derivative of eta III} to eliminate the $t$ derivatives of $\eta^{III}$, we can rewrite these terms as follows.
\begin{equation}\label{definition of phi III}
\begin{cases}
 \psi_1(y;\eta^{III})=
  \mathcal{C}_1(y)\partial_x^3\eta^{III}+\mathcal{C}_2(y)\delta\partial_x^4\eta^{III}+\cdots
+\mathcal{C}_7(y)\delta^6\partial_x^{9}\eta^{III}+N_1^{III},\\
 \psi_2(y;\eta^{III})=
\mathcal{C}_8(y)\partial_x^3\eta^{III}+\mathcal{C}_9(y)\delta\partial_x^4\eta^{III}+\cdots
+\mathcal{C}_{15}(y)\delta^7\partial_x^{10}\eta^{III}+N_2^{III},\\
 \phi_1(\eta^{III})=\mathcal{C}_{16}\partial_x^3\eta^{III}+\mathcal{C}_{17}\delta\partial_x^4\eta^{III}+\cdots
+\mathcal{C}_{21}\delta^5\partial_x^8\eta^{III}+N_3^{III},\\
 \phi_2(\eta^{III})=\mathcal{C}_{22}\partial_x^3\eta^{III}
+\mathcal{C}_{23}\delta\partial_x^4\eta^{III}+\cdots
+\mathcal{C}_{26}\delta^4\partial_x^7\eta^{III}+N_4^{III},\\
 \phi_3(\eta^{III})=\mathcal{C}_{27}\partial_x^4\eta^{III}+\mathcal{C}_{28}\delta\partial_x^5\eta^{III}+\cdots
+\mathcal{C}_{30}\delta^3\partial_x^7\eta^{III}
+N_5^{III},
\end{cases}
\end{equation}
where  $\mathcal{C}_1, \ldots, \mathcal{C}_{15}$ are polynomials in $y$, $\mathcal{C}_{16}, \ldots, \mathcal{C}_{30}$ are constants, and $N_1^{III}, \ldots, N_5^{III}$ are collections of the nonlinear terms  of the form
\begin{equation}\label{nonlinear terms of phi III}
\dfrac{1}{\delta^3}\Phi_0(\delta\eta^{III}, \delta^2\partial_x\eta^{III},  \ldots, \delta^5\partial_x^4\eta^{III};y)\Phi_0(\delta^2\partial_x\eta^{III}, \ldots, \delta^{10}\partial_x^9\eta^{III};y).
\end{equation}
Here we generally denote polynomials of $\bm{f}$ by the same symbol $\Phi=\Phi(\bm{f})$ 
and $\Phi_0$ is such a function satisfying $\Phi_0(\bm{0})=0$. 
We also use such a function $\Phi_0$ depending also on $y\in[0,1]$ and denote it by $\Phi_0(\bm{f};y)$, that is, 
$\Phi_0(\bm{0};y)\equiv 0$. 
Let $(\eta^\delta, u^\delta, v^\delta, p^\delta)$ be the solution of \eqref{NS-2}--\eqref{BC0-2} and we set $H^{III}:=\eta^\delta-\eta^{III}$, $\bm{U}^{III}:=(u^\delta-u^{III}, \delta(v^\delta-v^{III}))^{\rm T}$, $\tilde{\mathscr{E}}_m^{III}:=\tilde{\mathscr{E}}_m(H^{III}, \bm{U}^{III})$, and so on.
We prepare several lemmas to proceed the error estimate.
\begin{lem}\label{lem estimate of Nm2}
Under the same assumption as Proposition \ref{uniform estimate}, for 
any $\epsilon>0$ there exists a positive constant $C_\epsilon$ such that we have
\begin{equation}\label{estimate of Nm2}
\mathscr{N}_m^2(\bm{U}^{III})(t)\le
\epsilon \tilde{\mathscr{F}_m}(t)+C_\epsilon\delta^4\tilde{E}_m(t)\tilde{F}_{m+1}(t),
\end{equation}
where $\mathscr{N}_m^2$ is the collection of nonlinear terms defined by \eqref{definition of Nm2}.
\end{lem}
\medskip

\noindent
{\it Proof.}\quad By the explicit form of $\bm{f}$, $h_1$, and $h^{(3)}_2$ (see \eqref{definition of non-homogeneous term} and Section 2),  we can obtain the desired estimate in the same but more easier way as proving  [16, Lemmas 5.11 and 5.12].\qquad $\square$

\medskip

\begin{lem}\label{estimate of Nm3}
Under the same assumption as Proposition \ref{uniform estimate}, for any $\epsilon>0$ there exists a positive constant $C_\epsilon$ such that  we have
\[
\mathscr{N}_m^3(H;\eta^{III})(t)\le \epsilon\tilde{\mathscr{F}}_m(t)+C_\epsilon \delta^5|\eta_x^{III}(t)|_{m+12}^2, 
\]
where $\mathscr{N}_m^3$ is the collection of nonlinear terms defined by \eqref{definition of Nm3}.
\end{lem}
\noindent
{\it Proof.}\quad By  the well-known inequalities $\|\partial_x^k(fg)\|\lesssim\|f\|_{L^\infty}\|\partial_x^kg\|+\|g\|_{L^\infty}\|\partial_x^kf\|$ 
and $\|\partial_x^k\Phi_0(\bm{f};y)\|\le C(\|\bm{f}\|_{L^\infty})\|\partial_x^k\bm{f}\|$, \eqref{t-derivative of eta III}--\eqref{nonlinear terms of phi III} lead to $\sum_{k=0}^m\mathscr{N}_0^3(\partial_x^kZ_3)\lesssim \big(1+|\eta^{III}|_{m+12}^2\big)\delta^5|\eta_x^{III}|_{m+12}^2$.
Moreover, by Poincar\'e's inequality and \eqref{nonlinear terms of phi III}, we see that
$|(\partial_x^k H, \delta^3\partial_x^k\phi_3)_\Gamma|\le\epsilon\delta|\partial_x^kH_x|_0^2+C_\epsilon\delta^5|\partial_x^k\phi_3|_0^2 \le\epsilon\tilde{\mathscr{F}}_m+C_\epsilon \big(1+|\eta^{III}|_{m+12}^2\big)\delta^5|\eta_x^{III}|_{m+12}^2$.
These together with Lemma \ref{energy estimate app} imply the desired inequality.\qquad$\square$

\medskip


\begin{lem}\label{estimate of Nm1}
Under the same assumption as Proposition \ref{uniform estimate}, for any $\epsilon>0$ there exists a positive constant $C_\epsilon$ such that we have
\begin{align}
\mathscr{N}_m^1(H^{III}, \bm{U}^{III}, P^{III}; \eta^{III})(t)\le& (C_\epsilon\tilde{E}_2(t)+\epsilon)\tilde{\mathscr{F}}^{III}_m(t)\label{est of Nm1}+C_\epsilon\big\{\tilde{E}_m(t)\tilde{\mathscr{F}}^{III}_2(t)+\delta^4\tilde{E}_m(t)\tilde{F}_{m+1}(t)\\
&\ +\delta^5|\eta_x^{III}(t)|_{m+12}^2+(\tilde{F}_{m}(t)+\delta|\eta_x^{III}(t)|_{m+12}^2)\tilde{\mathscr{E}}_m^{III}(t)\big\},\notag
\end{align}
where $\mathscr{N}_m^1$ is the collection of nonlinear terms defined by \eqref{definition of Nm1}.
\end{lem}

\noindent
{\it Proof.}\quad In this proof, we omit the symbol $III$ appeared in a superscript of solutions for simplicity. By \eqref{definition of non-homogeneous term}, \eqref{tildeg}, and \eqref{definition of error}, we see that $\bm{F}$ is consist of terms of the form
\[
\begin{cases}
\delta\Phi_0(\eta^\delta, \delta\eta_x^\delta;y)(\nabla_\delta U_y, \nabla_\delta P)+\delta^2(\eta^\delta)^2(U_{yy}, P_y),\\
\delta\Phi_0(\eta^\delta, \delta\eta_x^\delta, u^\delta;y)(\delta V, \delta U_x),\\
\delta\Phi_0(\delta\eta_x^\delta, \delta\eta_t^\delta, \delta v^\delta;y)(U,  U_y),\\
\delta\Phi_0(\eta, \bm{u}, \nabla_\delta\bm{u}, \nabla_\delta u_{y}, \nabla_\delta p ;y)(\delta H_x, \delta H_t, U, \delta V),\\
\delta^2\eta^\delta(u_{yy}+p_y)H
\end{cases}
\]
and that
$G_2=\delta^2\{\eta^\delta(2H_x+U_x)+\eta_x^\delta U+(2\eta_x+u_x)H+uH_x\}$,  $G_3=\delta^2\{(\eta^\delta)^2H_x+(\eta^\delta+\eta)\eta_x H\}$, and $bH=2\delta\big(\delta(\eta^\delta)^2+(2+\delta\eta)\eta^\delta+\delta\eta^2+2\eta\big)H$.
Note that using \eqref{approximate solution III} and \eqref{t-derivative of eta III}, we can express the approximate solutions $\bm{u}, \nabla_\delta\bm{u}, u_{yy}$, and $\nabla_\delta p$ in terms of $\eta$ and its $x$ derivatives.  
In view of these,  by putting 
\[
\begin{cases}
\Phi^1=\Phi(\eta^\delta, \delta\eta_x^\delta, \delta\eta_t^\delta, \delta^2\eta_{xx}^\delta, \delta^2\eta_{tx}^\delta,\bm{u}^\delta ;y), \\
\Phi^2=\Phi(\delta\eta_x^\delta, \delta\eta_t^\delta, \delta^2\eta_{xx}^\delta, \delta^2\eta_{tx}^\delta, \delta^2\eta_{tt}^\delta, \delta v^\delta, \delta \bm{u}_x^\delta, \delta \bm{u}_t^\delta ;y),\\
\Phi^3=\Phi(\eta^\delta, \delta\eta_x^\delta;y),\\
\Phi^4=\Phi(\eta, \delta\eta_x, \ldots, \delta^{10}\partial_x^{10}\eta;y),
\end{cases}
\] 
\[
\begin{cases}
W:=(\delta H_x, \delta H_t, \delta^2 H_{xx}, \delta^2 H_{tx}, \delta^3 H_{xxx}, \delta V, \delta\bm{U}_x, \delta\bm{U}_t,  \delta\nabla_\delta U_x, \delta\nabla_\delta U_t, \nabla_\delta U_y, \nabla_\delta U_{xy},\\
\quad\qquad\nabla_\delta P, \nabla_\delta P_x, \delta U_x|_\Gamma, \delta U_t|_\Gamma, \delta^2 U_{xx}|_\Gamma, \delta^{5/2}|D_x|^{5/2} U|_\Gamma),\\
Q:=(H, \delta H_x, \delta H_t, \delta^2 H_{xx}, \delta^2 H_{tx}, \delta^3 H_{xxx}, \bm{U}, \nabla_\delta\bm{U}, \delta\bm{U}_t, U|_\Gamma),
\end{cases}
\]
it suffices to estimate
\[
\begin{cases}
I_1=\delta\|\partial_x^k(\Phi_0^1W)\|^2,\\
I_2=\delta\|\partial_x^k(\Phi_0^2Q)\|^2,\\
I_3=\delta^3|(\partial_x^k(\eta^\delta \bm{U}_{tx}), \partial_x^kV_t)_\Gamma|,\\
I_4=\delta^2|(\partial_x^k(\Phi_0^3\nabla_\delta U_{ty}), \partial_x^k\bm{U}_t)_\Omega|,\\
I_5=\delta^2|(\partial_x^k(\Phi_0^3\nabla_\delta P_t)_, \partial_x^k\bm{U}_t)_\Omega|,\\
I_6=\delta\|\partial_x^k(\Phi^1\Phi_0^4Q)\|^2,\\
I_7=\delta^4|(\partial_x^k(\Phi_0^4 H_{tt}), \partial_x^kV_t)_\Gamma|,\\
I_8=\delta^6{\rm W}|(\partial_x^kH_{xxxx}, \partial_x^kG_{3xx})_\Gamma|
\end{cases}
\]
for $0\le k\le m$.

By Proposition \ref{uniform estimate} and  $\|(u, v)\|_{L^{\infty}}\lesssim\|(u_y, v_y)\|+\|(u_{xy}, v_{xy})\|$ thanks to the boundary condition $u|_{y=0}=v|_{y=0}=0$, we obtain
\begin{align}
&\|\Phi_0^1\|^2_{L^\infty}\lesssim \tilde{E}_2,\quad \|\partial_x^k\Phi_0^1\|^2+\|\partial_x^k\Phi^1_{0y}\|^2\lesssim \tilde{E}_m,\label{phi1}\\
&\|\Phi_0^2\|^2_{L^\infty}\lesssim \tilde{F}_2,\quad \|\partial_x^k\Phi_0^2\|^2+\|\partial_x^k\Phi^2_{0y}\|^2\lesssim \tilde{F}_m\label{phi2}.
\end{align}
In the same way as the proof of Lemma \ref{estimate of Nm3}, we have
\begin{align}
\delta\|\Phi_0^4\|_{L^\infty}^2\lesssim\delta|\eta_x|_{m+12}^2,\quad\delta(\|\partial_x^k\Phi_0^4\|^2+\|\partial_x^k\Phi^4_{0y}\|^2)\lesssim\delta|\eta_x|_{m+12}^2,\quad |\Phi_0^4|_{m-\frac12}^2\lesssim|\eta|_{m+12}^2.\label{phi4}
\end{align}
On the other hand, it is easy to see that
\begin{align}
&\|W\|^2+\|W_{x}\|^2\lesssim\tilde{\mathscr{F}}_2,\quad \|\partial_x^kW\|^2\lesssim\tilde{\mathscr{F}}_m,\label{wi1}\\
&\|Q\|^2+\|Q_{x}\|^2\lesssim\tilde{\mathscr{E}}_2,\quad \|\partial_x^kQ\|^2\lesssim\tilde{\mathscr{E}}_m,\label{zi}
\end{align}
where we used the trace theorem $|f|_0^2+\delta||D_x|^{\frac12}f|_0^2\lesssim\|f\|^2+\delta^2\|f_x\|^2+\|f_y\|^2$ to estimate the term $\delta^5||D_x|^{\frac52}U|_0^2$. In the following, we often use the inequality
\begin{equation}\label{inequality2}
\|\partial_x^k(af)\|\lesssim\|a\|_{L^\infty}\|\partial_x^kf\|+(\|\partial_x^ka\|+\|\partial_x^k a_y\|)(\|f\|+\|f_x\|),
\end{equation}
which have been shown in [16, Lemma 5.2]. 

As for $I_1$, by \eqref{phi1}, \eqref{wi1}, and \eqref{inequality2}, we have $I_1\lesssim\tilde{E}_2\tilde{\mathscr{F}}_m+\tilde{E}_m\tilde{\mathscr{F}}_2$.
As for $I_2$, by \eqref{phi2}, \eqref{zi}, and \eqref{inequality2}, we have $I_2\lesssim\tilde{F}_m\tilde{\mathscr{E}}_m$.
As for $I_3$, by integration by parts, we have $I_3\lesssim C_\epsilon\delta^3|\eta^\delta U_{tx}|_{m-\frac12}^2+\epsilon\delta^3|V_t|_{m+\frac12}^2\le C_\epsilon(\tilde{E}_2\tilde{\mathscr{F}}_m+\tilde{E}_m\tilde{\mathscr{F}}_2)
+\epsilon\tilde{\mathscr{F}}_m$.
As for $I_4$, by integration by parts in $y$, we have 
\begin{align*}
I_4&\le C_\epsilon\delta^2\big(\|\partial_x^k(\Phi_0^3\nabla_\delta U_t)\|^2+\|\partial_x^k(\Phi_{0y}^3\nabla_\delta U_t)\|^2\big)\\
&\quad+\delta^3|(\partial_x^k(\Phi_0^3U_{tx}), \partial_x^k \bm{U}_t)_\Gamma|+\delta^2|(\partial_x^k(\Phi_0^3U_{ty}), \partial_x^k\bm{U}_t)_\Gamma|+\epsilon\delta\|\partial_x^k\bm{U}_{ty}\|^2\\
&\le I_{4,1}+I_{4,2}+I_{4,3}+\epsilon\tilde{\mathscr{F}}_m,
\end{align*}
where $I_{4,1}=C_\epsilon\delta^2\big(\|\partial_x^k(\Phi_0^3\nabla_\delta U_t)\|^2+\|\partial_x^k(\Phi_{0y}^3\nabla_\delta U_t)\|^2\big), I_{4,2}=\delta^3|(\partial_x^k(\Phi_0^3U_{tx}), \partial_x^k \bm{U}_t)_\Gamma|$, and $I_{4,3}=\delta^2|(\partial_x^k(\Phi_0^3U_{ty}), \partial_x^k\bm{U}_t)_\Gamma|$.
The estimates for $I_{4,1}$ and $I_{4,2}$ are reduced to the estimates for $I_1$ and $I_3$, respectively. Thus, taking into account that we can eliminate the term $U_y|_\Gamma$ in $I_{4,3}$ by the first equation in \eqref{BC1error}, this together with the estimates for $I_2$, $I_3$, $\delta^3h_1$, and $\delta^3\phi_1$ yields $I_4\le\epsilon\tilde{\mathscr{F}}_m+C_\epsilon\big\{\tilde{E}_2\tilde{\mathscr{F}}_m+\tilde{E}_{m}(\tilde{\mathscr{F}}_2+\delta^4\tilde{F}_{m+1}+|\eta|_{m+12}^2\delta^5|\eta_x|_{m+12}^2)\big\}$.
As for $I_5$,  it suffices to show the case of $k\ge1$ because we can treat easily the case of $k=0$. Integrating by parts in $x$, \eqref{phi1}, and \eqref{inequality2}, we have 
$I_5\le \epsilon\delta^3\|\partial_x^k\bm{U}_{tx}\|^2+C_\epsilon\delta\|\partial_x^{k-1}(\Phi^3_0\nabla_\delta P_t)\|^2\le\epsilon\tilde{\mathscr{F}}_m+C_\epsilon\big(\tilde{E}_2\tilde{\mathscr{F}}_m+\tilde{E}_m\tilde{\mathscr{F}}_2\big)$.
As for $I_6$, by \eqref{phi1},  \eqref{phi4}, \eqref{zi}, and \eqref{inequality2}, we have
\begin{align*}
I_6&\lesssim\delta\big\{\|\Phi_0^4\|_{L^\infty}^2(\|\partial_x^k\Phi^1\|^2+\|\partial_x^k\Phi^1_y\|^2)(\|Q\|^2+\|Q_{x}\|^2)\\
&\quad +\|\Phi^1\|_{L^\infty}^2(\|\partial_x^k\Phi_0^4\|^2+\|\partial_x^k\Phi^4_{0y}\|^2)(\|Q\|^2+\|Q_{x}\|^2)+\|\Phi^1\|_{L^\infty}^2\|\Phi_0^4\|_{L^\infty}^2\|\partial_x^kQ\|^2\big\}\\
&\lesssim (\tilde{E}_m+|\eta|_{m+12}^2)\delta|\eta_x|_{m+12}^2\tilde{\mathscr{E}}_m.
\end{align*}
As for $I_7$, it suffices to show the case of $k\ge1$ because we can treat easily the case of $k=0$. By the third equation in \eqref{BC1error}, integration by parts, and the trace theorem, we have
\begin{align*}
I_7&\le C_\epsilon\delta^4||D_x|^{\frac12}\partial_x^{k-1}(\Phi_0^4V_{t})|_0^2+C_\epsilon\delta^5|\partial_x^k(\Phi_0^4H_{xt}+\Phi_0^4G_{3t})|_0^2+C_\epsilon\delta^5|\delta^3\partial_x^k\phi_{3t}|_0^2\\
&\qquad+\epsilon\big(\delta^4||D_x|^{\frac12}\partial_x^kV_{t}|_0^2+\delta^3|\partial_x^kV_t|_0^2\big)\\
&\le I_{7,1}+I_{7,2}+I_{7,3}+\epsilon\tilde{\mathscr{F}}_m,
\end{align*}
where $I_{7,1}= C_\epsilon\delta^4||D_x|^{\frac12}\partial_x^{k-1}(\Phi_0^4V_{t})|_0^2$, $I_{7,2}=C_\epsilon\delta^5|\partial_x^k(\Phi_0^4H_{xt}+\Phi_0^4G_{3t})|_0^2$, and $I_{7,3}=C_\epsilon\delta^5|\delta^3\partial_x^k\phi_{3t}|_0^2$.
By the trace theorem, the second equation in \eqref{NSerror}, and \eqref{phi4}, we have
\begin{align*}
I_{7,1}&\lesssim|\Phi_0^4|_{m-\frac12}^2\delta^3|V_t|_{L^\infty}^2+\delta|\Phi_0^4|_{L^\infty}^2\delta^3||D_x|^{\frac12}\partial_x^{k-1}V_t|_0^2\\
&\lesssim|\Phi_0^4|_{m-\frac12}^2\delta^3\|U_{txx}\|^2+\delta|\Phi_0^4|_{L^\infty}^2(\delta^2\|\partial_x^{k}U_t\|^2+\delta^4\|\partial_x^kV_t\|^2)\\
&\lesssim|\eta|_{m+12}^2\tilde{\mathscr{F}}_2+\delta|\eta_x|_{m+12}^2\tilde{\mathscr{E}}_m.
\end{align*}
Recalling the explicit form of $G_3$, we see that the estimate of $I_{7,2}$ is reduced to $I_6$. Taking into account that we have already estimated $I_{7,3}$ in the proof of Lemma \ref{estimate of Nm3}, we obtain 
$I_7\le C_\epsilon\big\{|\eta|_{m+12}^2\tilde{\mathscr{F}}_2+(\tilde{E}_m+|\eta|_{m+12}^2)\delta|\eta_x|_{m+12}^2\tilde{\mathscr{E}}_m+|\eta|_{m+12}^2\delta^5|\eta_x|_{m+12}^2\big\}+\epsilon\tilde{\mathscr{F}}_m$.
As for $I_{8}$, integration by parts, \eqref{phi1}, and \eqref{phi4} lead to
\begin{align*}
\delta^6{\rm W}|(\partial_x^kH_{xxxx}, \partial_x^kG_{3xx})_\Gamma|&\le\epsilon(\delta^2{\rm W})^2\delta^2||D_x|^\frac72H|_m^2+C_\epsilon\delta^6||D_x|^\frac52G_3|_m^2\\
&\le \epsilon\tilde{\mathscr{F}}_m+C_\epsilon\big\{\delta^2(\tilde{F}_m+\delta|\eta_x|_{m+12}^2)\tilde{\mathscr{E}}_2+\tilde{E}_2\tilde{\mathscr{F}}_m\big\}.
\end{align*}

Therefore, by the boundedness of the terms $\tilde{E}_m$ and $|\eta|_{m+12}^2$ which comes from Proposition \ref{uniform estimate} and Lemma \ref{energy estimate app}, the proof is complete.\qquad$\square$

\bigskip

\begin{lem}\label{lem estimate for energy}
Under the same assumption as Proposition \ref{uniform estimate}, we have
\begin{align}
&\tilde{\mathscr{E}}_m^{III}(t)\lesssim \mathscr{E}_m^{III}(t)+\delta^4(\tilde{E}_{m+1}(t)+|\eta^{III}(t)|_{m+12}^2),\label{estimate for E}\\
& \tilde{\mathscr{F}}_m^{III}(t)\lesssim \mathscr{F}_m^{III}(t)+(\tilde{F}_{m}(t)+\delta|\eta_x^{III}(t)|_{m+12}^2)\tilde{\mathscr{E}}_m^{III}(t)\label{estimate for F}\\
&\phantom{\tilde{\mathscr{F}}_m^{III}(t)\lesssim}+\delta^4\tilde{E}_m(t)\tilde{F}_{m+1}(t)+\delta^5|\eta_x^{III}(t)|_{m+12}^2,\notag\\
&\mathscr{E}_m^{III}(t)\lesssim \mathscr{D}_m^{III}(t)+\delta^4.\label{estimate for D}
\end{align}
\end{lem}

\noindent
{\it Proof.}\quad  In view of  the discrepancy of non-homogeneous terms in the equations, modifying the proof of (6.2) in [16, Lemma 6.2], we obtain  \eqref{estimate for E}. 
Taking into account that we can eliminate $U_{yy}$ in $\mathscr{F}_m^{III}$ by using the first equation in \eqref{NSerror}, modifying the proof of  (6.3) in [16, Lemma 6.2],  it is not difficult to check that \eqref{estimate for F} holds. Moreover, modifying the proof of (6.10) in \cite{Ueno}, we obtain \eqref{estimate for D}.\qquad$\square$

\bigskip
\begin{lem}\label{estimate of initial energy}
Under the same assumption as Proposition \ref{uniform estimate}, we have
\[
\mathscr{D}_m^{III}(0)\lesssim \delta^4.
\]
\end{lem}
\noindent
{\bf Remark\ 5.1.}\ 
This lemma together with \eqref{estimate for D} yields
\begin{equation}\label{estimate of E0}
\mathscr{E}^{III}_m(0)\lesssim \delta^4.
\end{equation}

\noindent
{\it Proof.}\quad By the second and third equations in the compatibility conditions, we see that
\begin{align}
u_0(x,y)&=yu_{0y}(x,1)-\int_0^y\int_z^1 u_{0yy}(x,w){\rm d}w{\rm d}z\label{u0}\\
&=\big(2y\eta_0+4y\delta \eta_0^2+2y\delta^2\eta_0^3\big)+\delta y\big(-\delta v_{0x}+\delta^2h_1^{(0)}\big)-\int_0^y\int_z^1 u_{0yy}(x,w){\rm d}w{\rm d}z.\notag
\end{align}
It follows from \eqref{assumption III} and $\|(1+|D_x|)^{m+1}u_{yy}^{III}|_{t=0}\|\lesssim \delta$ (see the explicit form of $u^{III}$, that is, \eqref{definition of u0}--\eqref{pre approximate solution III} and \eqref{approximate solution III}) that $\|(1+|D_x|)^{m+1}u_{0yy}\|\lesssim\delta$.
Thus, by \eqref{u0}, the explicit form of $u^{III}$, \eqref{assumption III}, and the uniform estimate for $\delta^2|h_1^{(0)}|_{m+1}$ (see the proof of Lemma \ref{lem estimate of Nm2}), we obtain $\|(1+|D_x|)^{m+1}U|_{t=0}\|\lesssim\delta$. 
Combining this and the first equation in the compatibility conditions leads to $\|(1+|D_x|)^mV|_{t=0}\|\lesssim\delta$. Therefore, in view of the definition of $\mathscr{D}_m$ (see \eqref{definition of D}), using these and $H|_{t=0}=0$, we obtain the desired estimate.\qquad$\square$

%

\bigskip

\noindent
{\it Proof of \eqref{error estimate III} in Theorem \ref{main theorem}.}\quad
By Proposition \ref{uniform estimate}, Lemmas \ref{lem energy estimate NS}, \ref{lem estimate of Nm2}--\ref{estimate of Nm1}, and \eqref{estimate for E} and \eqref{estimate for F} in Lemma \ref{lem estimate for energy}, if $c_0$ and $\epsilon$ are sufficiently small, then we have
\begin{equation}\label{differential inequality}
\frac{{\rm d}}{{\rm d}t}\mathscr{E}^{III}_m(t)+\tilde{\mathscr{F}}^{III}_m(t)\le  C_1\big(\varphi_1(t)\mathscr{E}^{III}_m(t)+\tilde{E}_m(t)\tilde{\mathscr{F}}^{III}_2(t)+\delta^4\varphi_2(t)\big),
\end{equation}
where 
\begin{equation}\label{psi12}
\varphi_1(t)=\tilde{F}_{m}(t)+\delta |\eta_x^{III}(t)|_{m+12}^2, \quad \varphi_2(t)=\tilde{E}_{m}(t)\tilde{F}_{m+1}(t)+\delta|\eta_x^{III}(t)|_{m+12}^2.
\end{equation}
By considering the case of $m=2$ in \eqref{differential inequality} and using  Gronwall's inequality and Proposition \ref{uniform estimate}, if $c_0$ is sufficiently small, then we have 
$\mathscr{E}_2^{III}(t)+\int_0^t 
\tilde{\mathscr{F}}_2^{III}(s){\rm d}s\le\varphi_3(t)$, 
where 
\begin{equation}\label{psi3}
\varphi_3(t)=\mathscr{E}_2^{III}(0)\exp{\bigg(C_1\int_0^t \varphi_1(s){\rm d}s\bigg)}+C_1\int_0^t\delta^4\varphi_2(s)\exp{\bigg(C_1\int_s^t\varphi_1(\sigma){\rm d}\sigma\bigg)}{\rm d}s,
\end{equation}
which leads to
\begin{equation}\label{time integration of F2}
\int_0^t \tilde{\mathscr{F}}_2^{III}(s){\rm d}s\le\varphi_3(t).
\end{equation}
Note that by Proposition \ref{uniform estimate} and Lemma \ref{energy estimate app}, we have the exponential decay estimate for $\tilde{E}_{m+1}(t)$ and $|\eta^{III}(t)|_{m+13}^2$.
This together with \eqref{differential inequality}, Gronwall's inequality,  and $\delta\mathscr{E}_m^{III}\lesssim\tilde{\mathscr{F}}^{III}_m$ which comes from $|H|_0\lesssim|H_x|_0$ and $\|V\|\lesssim\|V_y\|=\|U_x\|$ (see \eqref{Em} and \eqref{Fm}) yields
\[
\mathscr{E}_m^{III}(t) \le\bigg\{\mathscr{E}_m^{III}(0)\exp{\bigg(C_1\int_0^t\varphi_1(s){\rm d}s\bigg)}+\varphi_4(t)\bigg\}{\rm e}^{-c\delta t},
\]
where
\begin{equation}\label{psi4}
\varphi_4(t)=C_1\int_0^t\big(\tilde{\mathscr{F}}_2^{III}(s)+\delta^4\tilde{F}_{m+1}(s)\big)\exp{\bigg(C_1\int_s^t\varphi_1(\sigma){\rm d}\sigma\bigg)}{\rm d}s.
\end{equation}
Combining  the above inequality and \eqref{estimate for E} and \eqref{estimate for D} in Lemma \ref{lem estimate for energy}, we obtain
\begin{equation}\label{integral inequality}
\tilde{\mathscr{E}}_m^{III}(t)\le C_2\big(\delta^4+\mathscr{D}_m^{III}(0)+\varphi_4(t)\big){\rm e}^{-c\delta t}.
\end{equation}
Here, recalling the definition $\eta^{III}(x, t)=\zeta^{III}(x-2t, \varepsilon t)$ and the assumption $\varepsilon=\delta$ and using  Lemma \ref{energy estimate app}, we have  $\int_0^\infty\delta|\eta^{III}_x(t)|_s^2{\rm d}t=\frac{1}{\varepsilon}\int_0^\infty\delta|\zeta^{III}_x(\tau)|_s^2{\rm d}\tau\lesssim|\eta_0|_s$. By this,  the integrability of $\tilde{F}_{m+1}$ which comes from Proposition \ref{uniform estimate}, and \eqref{estimate of E0}, we have $\varphi_3(t)\lesssim\delta^4$ (see \eqref{psi12} and \eqref{psi3}). This together with \eqref{time integration of F2} leads to $\varphi_4(t)\lesssim\delta^4$ (see \eqref{psi4}). Combining this, \eqref{integral inequality}, and Lemma \ref{estimate of initial energy}, we have
\begin{equation}\label{uniform estimate for error}
\tilde{\mathscr{E}}_m^{III}(t)\le C_3\delta^4{\rm e}^{-c\varepsilon t},
\end{equation}
which implies $\mathcal{D}(t; \zeta^{III}, u^{III}, v^{III}, p^{III})\lesssim\delta^4{\rm e}^{-c\varepsilon t}$ (see \eqref{difference} and \eqref{tildeEm}). Here, we used $\|V\|\lesssim\|V_y\|=\|U_x\|$.
Moreover, by taking into account the equality $P(x, y,t)=P(x,1, t)-\int_y^1 P_y(x, z, t){\rm d}z$ and using the second  equation in \eqref{NSerror}, the second equation in \eqref{BC1error}, and the uniform estimate \eqref{uniform estimate for error}, we easily obtain $\|(1+|D_x|)^{m}(p^\delta-p^{III})(t)\|^2\lesssim\delta^4{\rm e}^{-c\varepsilon t}$. Note that in the case of $O(\delta^{-1})\le{\rm W}\le O(\delta^{-2})$ we can estimate the term $\frac{\delta^2{\rm W}}{\sin\alpha}\partial_x^mH_{xx}$ which comes from the second equation in \eqref{BC1error} by $\tilde{\mathscr{E}}_{m+1}^{III}$.
Therefore,  the proof of \eqref{error estimate III} in Theorem \ref{main theorem} is complete.\qquad$\square$


\bigskip

We proceed to prove \eqref{error estimate I}, \eqref{error estimate II}, and \eqref{error estimate IV}. Let $\zeta^{I}$, $\zeta^{II}$, and $\zeta^{IV}$ be the solution for \eqref{zeta I}, \eqref{zeta II}, and \eqref{zeta IV}, respectively under the initial condition $\zeta^{I}|_{\tau=0}=\zeta^{II}|_{\tau=0}=\zeta^{IV}|_{\tau=0}=\eta_0$.
We put $\eta^{I}(x, t):=\zeta^{I}(x-2t, \varepsilon t)$, $\eta^{II}(x, t):=\zeta^{II}(x-2t, \varepsilon t)$, $\eta^{IV}(x, t):=\zeta^{IV}(x-2t, \varepsilon t)$ and 
\[
\begin{cases}
u^{I}(x,y,t):=u_0^{I}(y;\eta^{I}(x,t))+\delta u^{I}_1(y;\eta^{I}(x,t)),\\
v^{I}(x,y,t):=u^{I}_0(y;\eta^{I}(x,t))+\delta v^{I}_1(y;\eta^{I}(x,t)),\\
p^{I}(x,y,t):=p^{I}_0(y;\eta^{I}(x,t))+\delta p^{I}_1(y;\eta^{I}(x,t)),
\end{cases}
\]
\[
\begin{cases}
u^{II}(x,y,t):=u_0^{II}(y;\eta^{II}(x,t))+\delta u^{II}_1(y;\eta^{II}(x,t)),\\
v^{II}(x,y,t):=u^{II}_0(y;\eta^{II}(x,t))+\delta v^{II}_1(y;\eta^{II}(x,t)),\\
p^{II}(x,y,t):=p^{II}_0(y;\eta^{II}(x,t))+\delta p^{II}_1(y;\eta^{II}(x,t)),
\end{cases}
\]
\[
\begin{cases}
u^{IV}(x,y,t):=u_0^{IV}(y;\eta^{IV}(x,t))+\delta u^{IV}_1(y;\eta^{IV}(x,t))+\delta^2 u^{IV}_2(y;\eta^{IV}(x,t)),\\
v^{IV}(x,y,t):=u^{IV}_0(y;\eta^{IV}(x,t))+\delta v^{IV}_1(y;\eta^{IV}(x,t))+\delta^2 v^{IV}_2(y;\eta^{IV}(x,t)),\\
p^{IV}(x,y,t):=p^{IV}_0(y;\eta^{IV}(x,t))+\delta p^{IV}_1(y;\eta^{IV}(x,t))+\delta^2 p^{IV}_2(y;\eta^{IV}(x,t)),
\end{cases}
\]
where $u^{I}_0, v^{I}_0, p^{I}_0, \ldots$ were defined by \eqref{pre approximate solution I}--\eqref{pre approximate solution IV}. In view of this, by applying the same argument as showing \eqref{error estimate III}, it is not difficult to check that \eqref{error estimate I}, \eqref{error estimate II}, and \eqref{error estimate IV} holds. Therefore, the proof of Theorem \ref{main theorem} is complete.\qquad $\square$


\bigskip
Hiroki Ueno \par
Department of Mathematics, Faculty of Science and Technology, Keio University, \par
3-14-1 Hiyoshi, Kohoku-ku, Yokohama 223-8522, Japan. \par
E-mail: hueno@math.keio.ac.jp

\bigskip
Tatsuo Iguchi\par
Department of Mathematics, Faculty of Science and Technology, Keio University, \par
3-14-1 Hiyoshi, Kohoku-ku, Yokohama 223-8522, Japan. \par
E-mail: iguchi@math.keio.ac.jp

\end{document}